\documentclass[reqno,a4paper,11pt]{amsart}
\usepackage[utf8]{inputenc}

\usepackage{algorithm}
\usepackage{algpseudocode}
\usepackage[english]{babel}

\usepackage[utf8]{inputenc}
\usepackage[T1]{fontenc}
\usepackage[normalem]{ulem}

\usepackage{amssymb}
\usepackage{amsmath,amsfonts}
\usepackage[tikz]{bclogo}
\usetikzlibrary{arrows, arrows.meta}
\usepackage[top=3cm,bottom=3cm,left=2.5cm,right=2.5cm]{geometry}

\usepackage{lmodern}
\usepackage{textcomp}
\usepackage{mathtools, bm}
\usepackage{amssymb, bm}
\usepackage[unq]{unique}
\usepackage{enumerate}
\usepackage{comment}
\usepackage[breaklinks]{hyperref}

\usepackage[numbers]{natbib}  
\usepackage{cleveref}
\usepackage{graphicx} % Required for inserting images

\usepackage{amsthm}

\usepackage{tikz}
\usepackage{caption}
\usepackage{subcaption}
\newcommand{\abs}[1]{\left| #1 \right|}
\newcommand{\ept}{\text{ept}}

\newcommand{\bP}{\mathbb P}

\newtheorem{theorem}{Theorem}[section]
\newtheorem{lemma}[theorem]{Lemma}

\newtheorem{claim}[theorem]{Claim}
\newtheorem{corollary}[theorem]{Corollary}

\theoremstyle{definition}

\newtheorem{definition}[theorem]{Definition}
\Crefname{definition}{Definition}{Definitions} 
\Crefname{appendix}{Appendix}{Appendices}
\title{Tight bounds for expected propagation time of probabilistic zero forcing}
\author[Jelassi]{Mehdi Jelassi}

\address{Ecole Polytechnique Federale de Lausanne (EPFL), CH-1015 Lausanne,
Switzerland}
\email{mehdi.jelassi@epfl.ch}

\author[Portier]{Julien Portier}

\address{Ecole Polytechnique Federale de Lausanne (EPFL), CH-1015 Lausanne,
Switzerland}
\email{julien.portier@epfl.ch}

\author[Sarkar]{Rik Sarkar}

\address{Ecole Polytechnique Federale de Lausanne (EPFL), CH-1015 Lausanne,
Switzerland}
\email{rik.sarkar@epfl.ch}

\begin{document}

\begin{abstract}
We study the probabilistic zero forcing process, a probabilistic variant of the classical zero forcing process.
We show that for every connected graph $G$ on $n$ vertices, there exists an initial set consisting of a single vertex such that the expected propagation time is $n/2 + O(1)$. 
This result is tight and confirms a conjecture posed by Narayanan and Sun.
%This result is tight and confirms a conjecture posed by Narayanan and Sun [European Journal of Combinatorics 98 (2021)].
Additionally, we show tight bounds on the probabilistic throttling number, which captures the trade-off between the size of the initial set and the speed of propagation. Namely, we show that for every connected graph $G$ on $n$ vertices, there exists an initial set consisting of $O(\sqrt{n})$ vertices such that the expected propagation time is $O(\sqrt{n})$. 
This improves upon previous results by Geneson and Hogben, and confirms another conjecture by Narayanan and Sun.
\end{abstract}

\maketitle

%----------------------------------------------------------%
\section{Introduction}
\label{sec:intro}

%\mehdi{we should be consistent: sometimes we use "color" other times "colour". As we use the UK spelling, we should use colour everywhere} \julien{Except "color" and "neighbor", there should be no change in UK versus American spelling in this paper no? If yes, just choose one and make the necessary changes please.} \rik{I have used "color" and "neighbor"}
 Zero forcing is an iterative coloring process applied to a graph. 
It was originally introduced to address the maximum nullity problem in combinatorial matrix theory \cite{max-null2, max-null, max-null3, max-null4}, and has also found applications in quantum system control \cite{quantum}. 
Let $G$ be a graph with a vertex set $V$ and an edge set $E$. 
The zero forcing process on $G$ begins with an initial subset $S\subseteq V$ of blue-colored vertices, while the remaining vertices are colored white. 
A white vertex $v$ becomes blue if it is the only white neighbor of a blue vertex $u$: we then say that $u$ \emph{forces} $v$, denoted $u \to v$. 
If all vertices of $G$ eventually become blue through a finite number of repeated applications of this rule, then $S$ is called a \emph{zero forcing set}. 
The \emph{zero forcing number} $Z(G)$ of $G$ is defined as the smallest size of a zero forcing set of $G$.
Viewing zero forcing as a dynamical process on a graph, works by Chilakamarri et al. \cite{chilakamarri2012iteration}, Fallat et al. \cite{fallat2017zero}, and Hogben et al. \cite{hogben2012propagation} investigated the number of steps required for an initial subset of vertices $S$ to spread the color blue throughout the graph: this is known as the \emph{propagation time} of a zero forcing set $S$. \\
Probabilistic zero forcing, introduced by Kang and Yi \cite{kang2012probabilistic}, modifies the deterministic process by introducing randomness in the forcing step.
More precisely, given a current set \( B \) of blue vertices, each vertex \( u \in B \) attempts to force each of its white neighbors \( v \in \overline{B} \) to turn blue independently with probability
\[
\bP(u \to v) = \frac{|N[u] \cap B|}{\deg(u)},
\]
where \( N[u] \) is the closed neighborhood of \( u \). 
This rule is called the \textit{probabilistic color change rule}. 
\textit{Probabilistic zero forcing} is defined as the repeated application of this color change rule. \\
As noted in \cite{Hogben}, classical zero forcing has been used to model rumor spreading in social networks, however, due to the inherently unpredictable nature of human interactions, the probabilistic variant may offer a more realistic framework for such dynamics, thereby motivating its study.
Additionally, probabilistic zero forcing is closely related to the \emph{push} and \emph{pull} models, which have been well-studied in theoretical computer science as models for information diffusion in networks, see for instance \cite{daknama2021robustness,mehrabian2014randomized}. \\

A key parameter to investigate in this process is the \emph{expected propagation time}.
The \emph{propagation time} of a nonempty subset $S$ of vertices of $G$, denoted as $pt_{pzf}(G, S)$, is the random variable associated to the number of time steps in the probabilistic zero forcing process until all vertices in $G$ become blue, starting with the vertices in $S$ initially blue.
The \emph{expected propagation time} of $S$ for $G$ is the expected value of the propagation time of $S$, namely $\text{ept}(G,S)=\mathbb{E}[\text{pt}_{\text{pzf}}(G,S)]$. 
Another parameter of interest is the \emph{throttling number} $\text{th}_{\text{pzf}}(G)$ of a graph $G$, defined as 
\[
\text{th}_{\text{pzf}}(G) = \min_{S \subseteq V} \{ \left|S\right| + \text{ept}(G,S) \}.
\]
This parameter was initially introduced by Butler and Young \cite{Butler} as a way to measure the balance between resources allocated to accomplish a task (here the number $|S|$ of initially blue vertices) and time needed to accomplish the task ($\text{ept}(G,S)$ in our context).

Geneson and Hogben proved the following upper bound on the throttling number.
 \begin{theorem}[\cite{Hogben}, Theorem 6.5] 
 \label{theorem: bounds on th nb}
        For any connected graph $G$ on $n$ vertices, the probabilistic throttling number of $G$ satisfies $\text{th}_{\text{pzf}}(G) = O(\sqrt{n} \cdot \log^2 n)$.
\end{theorem}
They further note that this bound is optimal up to logarithmic factors, as the path $P_n$ on $n$ vertices satisfies $\text{th}_{\text{pzf}}(P_n) = \Omega(\sqrt{n})$. Narayanan and Sun \cite{Sun} conjectured that the $\log$ factors in \Cref{theorem: bounds on th nb} could be removed, i.e. that any connected graph $G$ on $n$ vertices satisfies $\text{th}_{\text{pzf}}(G) = O(\sqrt{n})$.
Our next result confirms this conjecture.

\begin{theorem}
\label{thm:Throttling}
    For any connected graph $G$ on $n$ vertices, the probabilistic throttling number of $G$ satisfies $\text{th}_{\text{pzf}}(G) = O(\sqrt{n})$.
\end{theorem}

More generally, our techniques extend to show that for every connected graph $G$ on $n$ vertices, and $s \geq 1$, there exists a set $S \subseteq V(G)$ of size at most $s$ such that $ept(G,S) = O(n/s)$.
However, for clarity and conciseness, we restrict our exposition to the proof of \Cref{thm:Throttling}. \\

In the second part of this paper, we investigate the expected propagation time when the initial set consists of a single vertex.
Given a graph $G$ with vertex set $V$, we define the \emph{expected propagation time} \( \text{ept}(G) \) for \( G \) as the minimum expected propagation time starting from a single vertex, namely
\[
\text{ept}(G) = \min_{v \in V} \mathbb{E}[\text{pt}_{\text{pzf}}(G,\{v\})].\]

Chan, Curl, Geneson, Hogben, Liu, Odegard, and Ross \cite{Chan} showed that $\ept(G) \leq \frac{e}{e-1}n$ for any connected graph $G$ on $n$ vertices.
This bound was improved by Narayanan and Sun to the following.

\begin{theorem}[\cite{Sun}, Theorem 4.11]
    For any connected graph $G$ on $n$ vertices, we have \( \ept(G) \le \frac{n}{2} + O(\log n) \).
\end{theorem}

Based on the suspicion that among all connected graphs on \( n \) vertices, the path \( P_n \) exhibits the slowest propagation, Narayanan and Sun further conjectured that their bound can be tightened to \( \ept(G) \le \frac{n}{2} + O(1) \).
We confirm their conjecture in the following result.

\begin{theorem} 
\label{ept improved}
    For any connected graph \( G \) on $n$ vertices, we have \( \ept(G) \le \frac{n}{2} + O(1) \).
\end{theorem}

We remark that this bound is tight, as for the path $P_n$ on $n$ vertices, we indeed have $\ept(P_n)=n/2+O(1)$ (see \cite{Hogben}). \\

The structure of the paper is as follows.
In \Cref{sec:Outline}, we present a brief outline of our proofs.
In Section~\ref{section preliminaries}, we introduce tools that underpin our later arguments.
In Section~\ref{section th nb}, we present the proof of \Cref{thm:Throttling}.
Finally, Section~\ref{section ept improved} is devoted to the proof of \Cref{ept improved}.

%----------------------------------------------------------%

\section{Proofs outline}
\label{sec:Outline}

Our proof of \Cref{thm:Throttling} relies on completely different ideas from the proof of the previous bound by Geneson and Hogben \cite{Hogben} in \Cref{theorem: bounds on th nb}.
One of the main challenges in proving \Cref{thm:Throttling} is that the probabilistic zero forcing process on a subgraph does not always relate to that of the host graph restricted to this subgraph.
Indeed, if the degrees of the vertices of the subgraph are not the same as the degrees in the host graph, then the probabilities of propagation may be different.
However, if we insist that the probabilistic zero forcing propagates in this subgraph only from the vertices which have the same degree in the subgraph as in the original graph, then one can relate the two processes.
Therefore, our approach consists of constructing a subgraph $T_v$ of $G$ for each $v$, with the following properties:
\begin{itemize}
    \item there exists a subset $D_v \subseteq V(T_v)$ such that $D_v$ is connected in $T_v$ and $N(D_v)=T_v \setminus D_v$,%every neighbor of every element of $D_v$ is in $V(T_v)$ and $D_v$ is a connected dominating set of $T_v$,
    \item $T_v$ is large but $D_v$ is not too big (namely $|T_v| \geq \sqrt{n}$ and $|D_v| \leq 2\sqrt{n}$), 
    %\mehdi{$D_v$ satisfies $\abs{D_v}>\sqrt{n}$ but $|D_v| \leq 2\sqrt{n}$}, namely $|D_v| \leq 2\sqrt{n}$, \mehdi{Shouldn't be $\abs{V(T_v)}>\sqrt{n}$ but $\le 2\sqrt{n}$ instead of $D_v$?}
    \item $v \in D_v$.
\end{itemize}
In this way, the degrees of vertices in $D_v$ are the same viewed in the graphs $T_v$ and $G$, and thus we may use them to show that, once a vertex in $T_v$ is blue, then all of $T_v$ becomes blue within $O(\sqrt{n})$ additional steps in expectation.
We then show that we can choose an initial set $S$ of blue vertices of size $O(\sqrt{n})$, such that each $T_v$ contains at least one blue vertex within $O(\sqrt{n})$ steps in expectation.
Overall, this shows that the graph becomes blue within $O(\sqrt{n})$ steps in expectation. \\

Our proof of \Cref{ept improved} builds on the work of Narayanan and Sun \cite{Sun}.
The main challenge in proving \Cref{ept improved} is that there exist some graphs $G$ for which the starting vertex must be chosen very carefully. 
A typical example is a path on $n$ vertices: in order to prove any bound of the form $n/2 + O(1)$, one must select a starting vertex near the centre of the path.
Narayanan and Sun observed that, at every step of probabilistic zero forcing on a graph $G$, the expected number of blue vertices increases by at least $2$, unless $G$ contains some specific structure, namely a cut-vertex or a \emph{cornerstone} (a particular cutset of size $2$, formally defined in the next section).
As a result, proving any bound close to $n/2$ requires careful treatment of these structures.
Taking these observations into account, Narayanan and Sun chose as a starting point a cut-vertex or cornerstone, whose removal separates the graph into two subgraphs with vertex sets $S$ and $T$, with no edges between them, and whose sizes are as close as possible, so that the process starts "near the centre of the graph".
Assuming that $|T| \geq |S|$, their proof then proceeds in three main steps:
\begin{itemize}
    \item First they show that in expected time $O(\log n)$, all neighbors of the starting cut-vertex or cornerstone become blue.
    \item Next, they show that in expected time $(|T|-|S|)/2$, at least $|T|-|S|$ vertices in $T$ become blue, due to the absence of cut-vertex or cornerstone during those steps by the careful choice of the starting point.
    \item Finally, they consider the simultaneous propagation of the probabilistic zero forcing through $S$ and the remaining white vertices in $T$, and show that in expectation both processes terminate after $|S|+O(1)$ steps.
\end{itemize}
This strategy gives an expected total number of steps of 
\begin{align*}
    O(\log n)+(|T|-|S|)/2+ |S|+O(1)= O(\log n)+(|T|+|S|)/2=n/2+O(\log n).
\end{align*}
We now explain how our approach refines this argument to obtain the tight bound in \Cref{ept improved}.
Our main idea consists of exploiting some of the $d$ neighbors of the starting cut-vertex or cornerstone that become blue during the initial $O(\log d)$ steps (in expectation) in Narayanan and Sun's algorithm.
If at least $d/2$ of these neighbors lie in $T$, this provides a "head-start" of $d/2$ vertices during the second phase of their strategy, directly giving the desired bound.
%\mehdi{we no longer need to distinguish the cases $d_T\ge \frac{d}{2}$ and $d_T\le \frac{d}{2}$. Therefore, we need to adapt this paragraph.} \julien{This is ok to keep it, as this is just to convey the main idea of the proof.}
On the other hand, if at least $d/2$ of these neighbors lie in $S$, this speeds up the propagation time in $S$ during the third phase of Narayanan and Sun's strategy.
However, since the bound of this phase depends on the maximum between the propagation time in $S$ and the propagation time in the remaining white vertices in $T$, this does not immediately yield an improved bound.
To overcome this issue, we propose a modified selection of the initial cut-vertex or cornerstone: one that optimises a carefully designed function reflecting the subtleties identified above.
Then following a strategy in the spirit of Narayanan and Sun and a careful analysis combining their methods with some new ideas, we show that this approach yields the tight bound $n/2+O(1)$ stated in \Cref{ept improved}.

\section{Preliminaries}
\label{section preliminaries}
Throughout this paper, we use classic graph theory notation.
For a graph $G=(V,E)$, and a subset \( S \subseteq V \), we denote by \( G[S] \) the induced subgraph of \( G \) on the vertex set \( S \), by \( \deg_S (v) \) the number of neighbors of \( v \) in \( S \) and by \( \deg(v) \) the degree of $v$.
Moreover, the \emph{open neighborhood} of $S$, denoted by $N(S)$, is the subset of vertices of $V \setminus S$ that have at least one neighbor in $S$. 
The \emph{closed neighborhood} $N[S]$ of $S$ is defined as $N[S] = N(S) \cup S$. 
A vertex \( v \) of a graph \( G \) is called a \emph{cut-vertex} if $G\setminus v$ is not connected. 
More generally, a subset $S$ of vertices of a graph \( G \) is called a \emph{cut-set} if $G\setminus S$ is not connected.\\

For $p \in [0,1]$, we let $Be(p)$ be a Bernoulli random variable of parameter $p$ (taking value $1$ with probability $p$ and $0$ with probability $1-p$).
We start by recalling the Chernoff bounds (see for instance Theorem 2.8 in \cite{JLR}).

\begin{theorem}
\label{thm:Chernoff}
    For each \( i \in [n] \), let \( X_i \sim \mathrm{Be}(p_i) \) for some \( p_i \in [0,1] \), and suppose that the random variables \( X_1, X_2, \ldots, X_n \) are mutually independent.  
    Let $X = \sum_{i=1}^n X_i$ and $\mu = \mathbb{E}[X]$.
    Then for every $t \geq 0$, we have
    \begin{align*}
        \mathbb{P}(X \geq \mu +t) \leq \exp\left(-\frac{t^2}{2(\mu+t/3)}\right),
    \end{align*}
    and
    \begin{align*}
        \mathbb{P}(X \leq \mu -t) \leq \exp\left(-\frac{t^2}{2\mu}\right).
    \end{align*}
\end{theorem}

The following result is an easy consequence of the Chernoff bounds.

\begin{lemma}
 \label{lem: high expectation after truncation}
 Suppose $X$ is a sum of independent Bernoulli random variables, satisfying $\mathbb{E}[X] \geq 2 + \frac{1}{5}$. Then we have $\mathbb{E}[\min(X,25)] \geq 2$.
\end{lemma}

\begin{proof}
    Let $\mathbb{E}[X]= \mu= 2 + \epsilon$, where $\epsilon \geq \frac{1}{5}$. If $\epsilon \geq 5$, then by \Cref{thm:Chernoff}, we have
    \begin{equation*}
        \mathbb{P}(X < 3) \leq e^{-\frac{(\epsilon -1)^{2}}{2(2+\epsilon)}}< \frac{1}{3}.
    \end{equation*}
    Therefore, we have $\mathbb{P}(X \geq 3)>\frac{2}{3}$, which in turn implies that $\mathbb{P}(\min(X,3) \geq 3) > \frac{2}{3}$. Thus, we have $\mathbb{E}[\min(X,25)] \geq \mathbb{E}[\min(X,3)] \geq 2$. Now suppose that $\frac{1}{5} \leq \epsilon <5$. Then, by the tail sum formula for the expectation, it follows that
    \begin{equation*}
       2+\epsilon=\mathbb{E}[X]= \mathbb{E}[\min(X,25)]+ \sum_{i=1}^{\infty} \mathbb{P}(X \geq 25+i).
    \end{equation*}
    By \Cref{thm:Chernoff}, we have $\mathbb{P}(X \geq 25+i) \leq e^{\frac{-(25+i-\mu)^2}{25+i+\mu}}\leq e^{-\frac{25+i-\mu}{2}}$, for $\mu <7$. Hence, we have
    \begin{equation*}
      \sum_{i=1}^{\infty} \mathbb{P}(X \geq 25+i) \leq \sum_{i=1}^{\infty}e^{-\frac{25+i-\mu}{2}} \leq \epsilon, \text{ for } \frac{1}{5} \leq \epsilon <5.
    \end{equation*}
    This implies that $\mathbb{E}[\min(X,25)] \geq 2$, which finishes the proof.
\end{proof}

% We use the classic asymptotic notation. 
% For functions $f=f(n)$ and $g=g(n)$ (with $g$ non-negative), we write $f=O(g)$ if there is a constant $C$ such that $|f(n)|\le Cg(n)$, and we write $f=\Omega(g)$ if there is a constant $c>0$ such that $f(n)\ge cg(n)$ for sufficiently large $n$.
% We write $f=\Theta(g)$ if $f=O(g)$ and $f=\Omega(g)$, and we write $f=o(g)$ or $g = \omega(f)$ if $f(n)/g(n)\to 0$ as $n\to\infty$.
%We give the background needed for \Cref{section th nb}.
We use the common definition of a submartingale. %\julien{Is it the right definition? \mehdi{yes} We want to use it with filtration right?}
%\julien{Is this definition formal? Shouldn't we say that $M_n$ belong to the filtration generated by $(X_i)_{0\leq i \leq n}$?}
    \begin{definition}[Submartingale]
A sequence of random variables $(M_n)_{n=0}^{\infty}$ with finite absolute means is called a \emph{submartingale} with respect to another sequence of random variables $(X_n)_{n=0}^{\infty}$ if for all \( n \), \( M_n \) is measurable with respect to the filtration generated by \( X_0, X_1, \dots, X_n \) and 
\[
\mathbb{E}[M_{n+1} \mid X_0, X_1, \dots, X_n] \geq M_n.
\]
% The sequence \( M_0, M_1, \dots \) is called a \textbf{supermartingale} if all conditions are the same, except 
% \[
% \mathbb{E}[M_{n+1} \mid X_0, X_1, \dots, X_n] \leq M_n.
% \]
\end{definition}

We recall Azuma's inequality for submartingales (see for instance \cite{Wormald1999a}).

\begin{theorem}\label{Azuma inq}
If $(X_k)_{k=0}^{\infty}$ is a submartingale and there exist constants $c_k > 0$ such that almost surely
\[
|X_k - X_{k-1}| \leq c_k \quad \text{ for all } k,
\]
then for every $t > 0$ and $N>0$, we have
\[
\mathbb{P}(X_N \leq X_0 - t) \leq \exp\left(-\frac{t^2}{2 \sum_{k=1}^{N} c_k^2} \right). %\tag{2.30}
\]
\end{theorem}

We now give the definition of a stopping time. 

\begin{definition}[Stopping Time]
Let $(X_n)_{n\ge 0}$ be a sequence of random variables and \( T \) be a random variable taking values in \( \{0,1,2,\dots\} \). $T$ is called a \emph{stopping time} with respect to $(X_n)_{n\ge 0}$ if for each \( n \), we have that \( \{T \leq n\} \) is a measurable function of \( X_0, X_1, \dots, X_n \). That is, 
\[
\{T \leq n\} \in \sigma(X_0, X_1, \dots, X_n) \quad \text{for all } n.
\]
\end{definition}

Intuitively, a random variable \( T \) is a stopping time if it is possible to determine whether \( T \leq n \) at time $n$. %In other words, the decision to stop cannot depend on future values of the random variables \( X_{n+1}, X_{n+2}, \dots \).
We adopt the formulation presented in \cite{textbook} for Doob's Optional Stopping Theorem.

\begin{theorem} [Doob's Optional Stopping Theorem]
\label{thm:Doob}
Suppose that \( M_n \) is a submartingale with respect to \( X_n \) and that \( T \) is a stopping time with respect to \( X_n \). Suppose that there exists a constant \( c > 0 \) such that \( |M_{n+1} - M_n| \leq c \) for all \( n \) and further assume that \( \mathbb{E}[T] < \infty \). Then,
\[
\mathbb{E}[M_T] \geq \mathbb{E}[M_0].
\]
% If \( M_n \) is a submartingale and all else is equal, then 
% \[
% \mathbb{E}[M_T] \geq \mathbb{E}[M_0].
% \]
% If \( M_n \) is a supermartingale and all else is equal, then 
% \[
% \mathbb{E}[M_T] \leq \mathbb{E}[M_0].
% \]
\end{theorem}

%finally I don't need this lemma (Mehdi 22 05)
%\begin{lemma}[~\cite{Sun}, Lemma 3.5] \label{lemma 3.5 de N and S}
%Fix a blue vertex $v$ of $G$. Let $w \neq v$ be some other vertex such that the shortest distance between $v$ and $w$ is $s$. Then there exist explicit constants $C, C' > 0$ and $0 < \log(\alpha) < 1$ such that after $t + C' s + C s \log n$ steps, $w$ will be blue with probability at least $1 - \log(\alpha)^t$.
%\end{lemma}

%    Now, we give the background needed for \Cref{section ept improved}.

It will sometimes be useful in our proofs to consider the probability that a given subset of vertices is blue at a specific time step $t$. 
Given a graph $G$ with vertex set $V$, and a subset \( U \subseteq V \), we define \( P^{(t)}(G, S, U) \) as the probability that all vertices in \( U \) are blue after \( t \) steps of probabilistic zero forcing on $G$ starting with an initial set $S$ of blue vertices.
Moreover, define $P^{(t)}(G, S) = P^{(t)}(G, S, V)$ to be the probability that all vertices of \( G \) are blue  after \( t \) steps of probabilistic zero forcing on $G$ starting with initial set $S$ of blue vertices. 
Define \( \text{ept}(G, S, U) \) as the expected number of steps needed until all vertices in \( U \) are blue, noting that 
\[
\text{ept}(G, S) = \text{ept}(G, S, V).
\]
%\end{definition}

The following result from \cite{Hogben} claims that starting with a larger initial set of blue vertices cannot increase the expected propagation time of probabilistic zero forcing.
%In , the authors proved the following result which says that with a larger initial set of blue vertices, it is more likely to color the whole graph $G$. Hence, starting with a smaller initial set of blue vertices increases the expected propagation time of probabilistic zero forcing.
\begin{lemma}[{\cite[Proposition 4.1]{Hogben}}]\label{lemma 2.10 for Sun}
Suppose that $S \subseteq T$. Then, $P^{(\ell)}(G, S) \leq P^{(\ell)}(G, T)$ for every integer $\ell$, and thus $\operatorname{ept}(G, S) \geq \operatorname{ept}(G, T)$.
\end{lemma}
For our purposes, we will also need the following stronger result.
%Commentaire UTILE
%ok c'est stronger car c'est pas juste réduire la taille du initial set ob blue vertices mais c'est propager from each vertex less likley. (donc lemme 3.9 est un cas particulier où on garder les mêmes proba partout sauf à some vertices de T où on met à 0 pour réduire T à S qui est smaller.)
\begin{lemma}[{\cite[Lemma 2.11]{Sun}}]\label{lemma 2.11 for Sun}
Suppose that initially, some set $S \subseteq V$ is blue. Say that we follow some modified probabilistic process where at the $t^{th}$ step, $P_t(u \to v)$, the probability that $u$ converts $v$ to become blue at step $t$, is some function of $G, u, v,$ and $B_{t-1}$, the set of blue vertices after the $(t-1)^{th}$ step. 

In addition, suppose that
\[
P_t(u \to v) \leq \frac{|N[u] \cap B_{t-1}|}{\deg(u)}
\]
for all blue vertices $u$ and white neighbors $v$ of $u$, and that conditioned on $B_{t-1}$, the set of events $u \to v$ at step $t$ is independent. 

Then, for any $T \subseteq V$ and any $\ell \geq 1$, the probability that all vertices in $T$ are blue after time step $\ell$ is at most $P^{(\ell)}(G, S, T)$, the probability that all vertices in $T$ would be blue if we followed the normal probabilistic zero forcing process. 

Consequently, the expected amount of time until all vertices in $T$ are blue is at least $\operatorname{ept}(G, S, T)$.
\end{lemma}

It will also be useful in our proofs to state the following slightly more general version of \cite[Lemma 3.4]{Sun}.
The proof is essentially identical, and therefore is deferred to \Cref{Appendix-Bounds-Star}.

% Now, we show a slight variant of \Cref{color the star}. \julien{Make better transition}
% Given a blue vertex $v\in V$ of degree $d\ge 2$, as for any $t>\abs{N[v]}$, the following lemma is necessary as directly applying Lemma \ref{lem:concentration of each T_v} to $G[N[v]]$ only provides an exponential decay for the event $\text{Step}(G[N[v]],\{v\})\ge 2\abs{N[v]}+t$ whereas we aim to obtain an exponential decay for the event $\text{Step}(G[N[v]],\{v\})\ge 2\log d+t$.
%\julien{Write the next lemma for the stopping number of reaching any number of blue leaves. And then maybe write this as a simple corollary}

\begin{lemma}\label{lem: coloring d neighbours in a star}
There exist constants $C>0$ and $0<\alpha<1$ such that the following holds.
Let $ 2 \leq d \leq n$ be two integers. 
Let $H$ be a star graph with $n$ leaves, with its center vertex colored blue and all other vertices colored white. Then, in an instance of probabilistic zero forcing on $H$, after $t$ steps, the number of blue leaves is at least $d$ with probability at least $1-\alpha^t$, provided $t>C \log d$.
\end{lemma}

We obtain the following simple Corollary from Lemma \ref{lem: coloring d neighbours in a star}, whose proof can also be found in Appendix \ref{Appendix-Bounds-Star}.

\begin{corollary}\label{cor: color d neighbours in log d steps}
Let $ 2 \leq d \leq n$ be two integers. 
Let $H$ be a star graph with $n$ leaves. 
Let $\Gamma$ be the random variable associated with the number of time steps it takes for at least $d$ leaves to be forced blue in an instance of probabilistic zero forcing on $H$, with only the center vertex initially colored blue. Then we have 
\begin{align*}
    \mathbb{E}[\Gamma]=O(\log d).
\end{align*}
\end{corollary}

Finally we show the following simple lemma about set systems.

\begin{lemma}
 \label{lem:Set-System}
 Let $n\ge 4$ and $T_1,\dots, T_n$ be subsets of $[n]$ each of size at least $\sqrt{n}$. Then there exists a subset $S$ of $[n]$ of size $\lfloor2\sqrt{n}\rfloor$ such that $\cup_{i\in S}T_i$ intersects $T_j$ for each $j\in [n]$.
        
\end{lemma}
    \begin{proof}
    We may first suppose that each $T_i$ has size exactly $\left\lfloor \sqrt{n}\right\rfloor $.
   Let $S = \{a_1, \dots, a_{\left \lfloor 2\sqrt{n}\right \rfloor}\}$ be a set of $\lfloor 2\sqrt{n} \rfloor$ indices maximizing $\abs{\bigcup_{s \in S'} T_s}$ over all subsets $S' \subseteq [n]$ of size $\lfloor 2\sqrt{n} \rfloor$.
   It is enough to show that $\cup_{i\in S}T_i$ intersects each $T_j$ for $j\notin S$.
    Suppose by contradiction that $T_l\cap \bigg(\bigcup_{i\in S}T_i\bigg) = \emptyset$ for some $l\notin S$. 
    If there exists $j\in S$ such that $T_{j} \cap \bigg( \bigcup_{i\in S \setminus \{j \}} T_{i} \bigg) \neq \emptyset$, then we have $\abs{\bigcup_{s \in S \setminus \{j\}} T_s} > \abs{\bigcup_{s \in S} T_s} - \abs{T_j}$.
    %\mehdi{only the following is useful : $\abs{\bigcup_{s \in S \setminus \{j\}} T_s} > \abs{\bigcup_{s \in S} T_s} - \abs{T_j}$, without $\sqrt n$}
    Therefore, if we let $S' = S \setminus \{ j \} \cup \{ l \}$, then we have 
    \begin{align*}
        \abs{\bigcup_{s \in S'} T_s} > \abs{\bigcup_{s \in S} T_s} - \abs{T_j} + \abs{T_l} = \abs{\bigcup_{s \in S} T_s},
    \end{align*}
    which contradicts the definition of $S$.
    Therefore, for every $j \in S$, we have $T_{j} \cap \bigg( \bigcup_{i\in S \setminus \{j\}} T_{i} \bigg) = \emptyset$.
    Thus, all $T_{i}$ must be disjoint: $T_{i}\cap T_{j}=\emptyset$ for all distinct $i,j \in S$.
    But then
    \begin{align*}
        n < \lfloor 2\sqrt{n} \rfloor \cdot \lfloor \sqrt{n} \rfloor = \sum_{i \in S} \abs{T_{i}}=\abs{\bigsqcup_{\substack{i\in S}} T_{i} }\leq n,
    \end{align*} 
    which is a contradiction.
    \end{proof}

%\julien{We should probably remove those lemmas eventually, but let's keep them for now until I make up my mind.}
%The following lemmas bound the expected time of some steps of the algorithm proposed by Narayanan and Sun in \cite{Sun}.
%\begin{lemma}[~\cite{Sun}, Lemma 4.8]\label{star lemma}
%    Step 4 takes $O(\log \text{}n)$ time in expectation.
%\end{lemma}

%\begin{lemma}[~\cite{Sun}, Lemma 4.9]\label{twice speed lemma}
%    Step 5 takes at most $\frac{1}{2}(\left|T\right| - \left|S\right|)$ steps in expectation.
%\end{lemma}
%\begin{lemma}[~\cite{Sun}, Lemma 4.10]\label{attack the S sides lemma}
%    Step 6 takes $\left|S\right| + O(1)$ steps in expectation.
%\end{lemma}

\section{Upper bound for the Throttling number} \label{section th nb}

We start with a few additional definitions.
Given a graph $G$ with vertex set $V$, and two subsets $U, S \subseteq V$, we let $pt_{pzf}(G, S,U)$ be the random variable associated to the number of time steps required in the probabilistic zero forcing process until all vertices in $U$ become blue, starting with the vertices in $S$ initially blue. \\
As explained in \Cref{sec:Outline}, we will construct, for each vertex $v$, some subgraph containing $v$ on which the probabilistic zero forcing process will be easier to analyse.
We also need here some additional definitions.
Let $H$ be a connected subgraph of $G$.
A subset of vertices $C \subseteq V(H)$ is said to be a \emph{core} of $H$ if $H[C]$ is connected and $N_G(C)=V(H)\setminus C$. 

% \julien{Change all the sentences to make it more smooth}
% To prove \Cref{thm:Throttling}, we construct recursively a decomposition of a connected graph $G$ into connected subgraphs so that the probabilistic zero forcing process on each subgraph does relate to that of the host graph restricted to each subgraph. This is made precise by the definition of a core probabilistic zero forcing. 
%\rik{I have edited the definition below, please check}
    %We define leaves a subgraph $L$ to be  $L:=H\setminus C\subseteq H$ 
    %and we say $H$ is decomposed by its core $C$ and its leaves $L$ if $H=C\sqcup L$.

    \begin{definition}
    \label{def: v good}
        Given a graph $G$ on $n$ vertices and a vertex $v \in V(G)$, we say that a pair $(H,C)$ is \emph{$v$-good} if $H$ is a subgraph of $G$ and
        \begin{itemize}
            \item $C$ is a core of $H$,
            \item $v \in C$,
            \item $2\sqrt{n}\le\abs{V(H)}$,
            \item there exists a vertex $w \in C$ and a subset $W \subseteq N(w)$ such that $|V(H) \setminus W| < 2\sqrt{n}$ and $C \setminus \{w\}$ is a core of $G[V(H) \setminus (W \cup \{w\})]$ in $G \setminus \{w\}$. 
            % \mehdi{why not only $|V(H) \setminus W| < 2\sqrt{n}?$}
            % \julien{This is correct, so I am leaving this comment so that we make sure when re-reading it that this is taken into account}
            % \rik{I have taken this into account and modified the proof of the next lemma accordingly}
            %\rik{$C\setminus \{w\}$ is a core of $H \setminus W$?}
        \end{itemize}
    \end{definition}

The following result shows that such $v$-good pairs can be constructed for every vertex $v$ in a connected graph $G$. %provides the construction of a core $C$ of a connected graph with some additional desired properties.
\begin{lemma}\label{decomposition core leaves}
Let $G$ be a connected graph on $n\ge 4$ vertices.
For every $v\in V(G)$, there exists a $v$-good pair of $G$.
\end{lemma}

\begin{proof}
    We recursively construct vertex sets $C_{i}$ and induced subgraphs $H_i$ of $G$, such that $C_i$ is a core of $H_i$ and $v \in C_i$. Let $C_0 = \{v\}$ and $H_0=G[C_{0} \cup N(C_0)]$. Clearly, $C_0$ is a core of $H_{0}$ and $v \in C_{0}$. Now, given $C_{i}$ and $H_i$, we construct $C_{i+1}$ and $H_{i+1}$ as follows. If $C_i=V(G)$, then we stop. Otherwise, since $G$ is connected, $N(C_i) \neq \emptyset$. Pick an arbitrary vertex $w_{i+1} \in N(C_i)$. Let $C_{i+1}=C_i \cup \{w_{i+1}\}$ and $H_{i+1}=G[C_{i+1}\cup N(C_{i+1})]$. Clearly, $v \in C_i \subseteq C_{i+1}$. Since $H_i[C_i]$ is connected and $w_{i+1} \in N(C_i)$, then $H_{i+1}[C_{i+1}]$ is connected. Furthermore, $N(C_{i+1})=V(H_{i+1})\setminus C_{i+1}$ by construction, so $C_{i+1}$ is a core of $H_{i+1}$. \\
    Let $\ell$ be the smallest integer such that $|V(H_{\ell})| \geq 2 \sqrt n$. Such an $\ell$ must exist, since $|V(H_i)| \geq |C_{i}| = i+1$ for every $i$, and $n \geq 2 \sqrt n$ for $n \geq 4$. 
    We claim that $(C_{\ell}, H_{\ell})$ is a $v$-good pair. Clearly, $v \in C_{\ell}$, $C_{\ell}$ is a core of $H_{\ell}$ and $|V(H_{\ell})| \geq 2 \sqrt n$. If $\ell=0$, it is easy to see that $w=v$ and $W=N(v)$ satisfy the last condition in Definition \ref{def: v good}. Otherwise, if $\ell \geq 1$, we claim that $w=w_{\ell}$ and $W=V(H_{\ell})\setminus V(H_{\ell-1})$ satisfy the condition prescribed in Definition \ref{def: v good}. Note that $V(H_{\ell})\setminus V(H_{\ell-1})= N(C_{\ell}) \setminus N(C_{\ell-1})=N(w_{\ell}) \cap \overline{V(H_{\ell-1})} \subseteq N(w_{\ell})$. By our choice of $\ell$, $|V(H_\ell)\setminus W|=|V(H_{\ell -1})| < 2 \sqrt n$, and $C_{\ell -1}=C_{\ell} \setminus \{w_{\ell}\}$ is a core of $G[V(H_\ell) \setminus W ]=H_{\ell-1}$ in $G$. Since $w_{\ell} \in N(C_{\ell -1})$, $C_{\ell -1}$ is a core of $G[V(H_{\ell-1}) \setminus \{w_{\ell}\}]$ in $G \setminus \{w_{\ell}\}$.
\end{proof}

The next result bounds the tail of the propagation time of probabilistic zero forcing in a subgraph $H$ which has a core $C$, starting from a vertex in $C$.
% Intuitively, the following result asserts that the random variable which denotes the number of steps needed until $T_v$ is completely blue given it contains some blue vertex has a flat tail. \julien{Rewrite this sentence}
%\rik{I am trying to do this for arbitrary $v\in V(H)$, not necessaily in the core}
\begin{lemma}\label{lem:concentration of each T_v}
Let $G$ be a connected graph, and let $H$ be a connected subgraph of $G$ such that $C$ is a core of $H$.
Suppose $H$ has $h\ge 2$ vertices, and let $v\in C$. 
Let $\Gamma=pt_{pzf}(G,\{v\},V(H))$ be the propagation time of the probabilistic zero forcing process on $G$, starting from the initial blue set $\{v\}$, until $V(H)$ becomes completely blue.
Then, for any $t \geq h$, we have 
\[\mathbb{P}(\Gamma \geq 2h+t)\leq \exp \left( -\frac{t}{100} \right).\] 
%\mehdi{we rather have $\exp \left( -\frac{t}{24} \right)$}
% \julien{This is cleaner to keep round numbers, especially since we do not care about the multiplicative constant}
\end{lemma}
\begin{proof}
%Let $\Gamma_1$ denote the random variable associated with the number of steps taken until a vertex $u \in C$ turns blue, and let $\Gamma_2 = \Gamma - \Gamma_1$. If $v \in C$, then $\mathbb{P}(\Gamma_1 \geq t/2)=0$. On the other hand, if $v \in V(H) \setminus C$, then $v$ has a neighborin $C$. Then $\mathbb{P}(\Gamma_1 \geq t/2 )$ is at most $\mathbb{P}(pt_{pzf}(S_a,\{w_a\}) \geq t/2)$, where $S_a$ is a star on $d(v)+1$ vertices centered at $w_a$. We have $t/2 \geq \sqrt n =\omega(\log n)$, so for $n$ sufficiently large, $\mathbb{P}(pt_{pzf}(S_a,\{w_a\}) \geq t/2) \leq \alpha^{\frac{t}{2}}$, by Lemma \ref{lem: coloring d neighbours in a star}. After $\Gamma_1$ steps, a vertex $u \in C$ is colored blue, so from here on, 

We consider the alternative probabilistic process where only vertices in $C$ are allowed to propagate their color (still independently and with same probabilities as in the probabilistic zero forcing process on $G$). Let $\Gamma'$ denote the number of steps taken by this alternative process, starting with the vertex $v \in C$ initially blue, until all vertices in $H$ become blue. Then, by Lemma \ref{lemma 2.11 for Sun}, $\mathbb{P}(\Gamma \geq 2h+t) \leq \mathbb{P}(\Gamma' \geq 2h+t)$, and we focus on bounding the tail of $\Gamma'$.\\

Let $S_i$ be the set of blue vertices in $H$ at the $i$th step of the alternative process.
Let $(X_i)_{i=1}^{\infty}$ be the sequence of random variables defined in the following way.
%\julien{Change $v_j$ to $v^i$}\rik{done}
Let $v^i$ be the vertex of minimum index in $C$ such that $v^i \in S_i$ and $N(v^i) \not\subseteq S_i$
(in other words, at step $i$, $v^i$ is blue and has at least one white neighbor). Then we set $X_{i+1}=1$ if $v^i$ gave its color to at least one of its white neighbors at step $i+1$, and $X_{i+1}=0$ otherwise.
If there exists no such $v^i$, we also set $X_{i+1}=0$. 
We remark that if at step $i$, there exists no such $v^i$ as above, then $V(H)$ is already colored blue.
Indeed, if $V(H)$ was not already colored blue, then there must exist a white vertex $y \in V(H)$, and thus, as $C$ is a core of $H$ and $H$ is connected, there must exist a path $vv_1\dots v_{\ell}y$ such that all the $v_i$ are in $C$. Setting $v_0=v$ and $v_{\ell+1}=y$, and since $v$ is blue and $y$ is white, it follows that there exists $0 \leq k \leq \ell$ such that $v_k$ is blue and $v_{k+1}$ is white, as desired. %\mehdi{but to conclude as such, we need to assume that we propagate in the order of the labelling}\\
% \begin{Claim}
%     If $S_i \subseteq V(G)$ is such that , then $S_i=V(G)$.
% \end{Claim}

% \begin{proof}
    
% \end{proof}

% as for any $i$, there exist some vertices $v_j$ in $V(C)$ that are blue at step $i$ for $j\in \{1,\dots,\abs{S_i\cap V(C)}\}$, such that $X_i=1$ if $v_j$ gave its color to at least one of its white neighbors at step $i+1$. If there exists no such $v_j$, then set $X_i=0$. In particular, if $G$ is blue at step $i$ then there exists no such $v_j$ so that $X_i=0.$
% \\
% We explain why, if there exists no $v_j\in S_{i}\cap V(C)$ that has at least one white neighbor, then $G$ must be blue at step $i+1$. 

%If $G$ is not entirely colored in blue at step $i+1$, then there exists a blue vertex that is adjacent to at least one white vertex (since $G$ is connected) to which it has not yet propagated its color. But by definition of core probabilistic zero forcing, this blue vertex must be in $V(C)$.

% Indeed, no core vertex $v_j$ in $S_{i}$ has at least one white neighbor means that all neighbors of $v_j$, for each $j$, that are in the core $C$ are blue (here we use that $C$ is connected) and that
% $N(V(C))=V(T_v)\setminus V(C)$ is a set of blue vertices (here we use that $G$ is connected). Hence $G$ is blue. 
% %All paths from $C$ to $L$ are entirely colored in blue, and all neighbors of vertices in the subpath in $C$ as well (here we use that $C$ is connected). Hence $G$ is blue. 
% \\

We now condition on $S_0, \dots, S_i$, and let $v^i$ be the vertex of minimum index in $C$ such that $v^i \in S_i$ and $N(v^i) \not\subseteq S_i$.
Let $\{w_1, \dots, w_a\} = N(v^i) \cap \bar{S_i}$ (i.e. the set of white neighbors of $v^i$ after step $i$), and let $N=|N(v^i)|$.
Then $v^i$ propagates its color to each of $w_1, \dots, w_a$ independently with probability $\frac{1+N-a}{N}$.
Therefore,
\begin{align*}
    \bP(X_{i+1} = 0 | S_0, S_1, \dots, S_i) = \left(1-\frac{1+N-a}{N}\right)^a \leq \exp\left(-\frac{a(1+N-a)}{N}\right) \leq e^{-1},
\end{align*}
from which it follows that 
\begin{align}
\label{eq:Bound-Proba-New-Vertex-Blue}
    \bP(X_{i+1} = 1 | S_0, \dots, S_i) \geq 1-e^{-1}.
\end{align}
% Let $j\in \{1,\dots, \abs{S_{i}\cap V(C)}\}$ and $D_{i,j}$ be the event no blue vertex $v_j$ at step $i$ in $V(C)$ has given its color to one of its $h_{W,j}$ white neighbors $w_1,\dots,w_{h_{W,j}}$ at step $i+1$.

% Let $D_i$ be the event 
% no blue vertex $v_j$ at step $i$ in $V(C)$ has given its color to one of its white neighbors at step $i+1$ so that $\bigcap_{j=1}^{\abs{S_{i}\cap V(C)}} D_{i,j}=D_i$.

% As $D_i$ occurs if and only if $X_i=0$ and by independence, we have if $G$ is not entirely blue at step $i$ the following 
% \begin{align*}
% \bP(X_{i+1} = 0 | S_1, \dots, S_i) = \bP(D_i) 
% &= \bP\left( \bigcap_{j=1}^{\abs{S_{i} \cap V(C)}} D_{i,j} \right) \\
% &\leq \bP(D_{i,p}) = \left(1 - \frac{c_{i,p} + 1}{h_{W,p} + c_{i,p}} \right)^{h_{W,p}} \\
% &\leq \exp\left(-\frac{h_{W,p}(c_{i,p} + 1)}{h_{W,p} + c_{i,p}} \right) \\
% &\leq \exp(-1)
% \end{align*}

% where $c_{i,p}\geq 1$ and $h_{W,p}$ are the number of blue neighbors and white neighbors of $v_p$ respectively and the independence follows by definition of core probabilistic zero forcing on $G$; and where $v_p$ is some vertex in $S_{i}\cap V(C)$ that has at least one white neighbor i.e. $h_{W,p}\ge 1$ (as long as $G$ is not entirely colored in blue at step $i$, $v_p$ exists as explained in the first paragraph).
%(\mehdi{since some $v_j$'s may not have white neighbors but as long as $G$ is not blue, $v_p$ exists})

We now define the sequence of random variables $(M_i)_{i=0}^{\infty}$ such that $M_0 =0$, and $M_{i+1}=M_i+X_{i+1}-(1-e^{-1})$ if there is some white vertex in $H$ at step $i$ and $M_{i+1}=M_i$ otherwise. 
%Also, consider the following stopping time: the first time $i$ for which after $i$ steps, $G$ is colored entirely in blue.% If at step $i$, $G$ is entirely colored in blue, then $M_i=M_{i+1}$. 
If $S_i \neq V(H)$, from \eqref{eq:Bound-Proba-New-Vertex-Blue}, we have that 
\begin{align*}
\mathbb{E}[M_{i+1} \mid S_0, \dots, S_i] 
= M_i + \mathbb{E}[X_{i+1} - (1-e^{-1})\mid S_0, \dots, S_i] \ge M_i.
\end{align*}
Otherwise, if $S_i = V(H)$, then $\mathbb{E}[M_{i+1} \mid S_0, \dots, S_i] =M_i$. Thus, $(M_i)_{i=0}^{\infty}$ is a sub-martingale with respect to $(S_i)_{i=0}^{\infty}$. Moreover, it is easy to see that almost surely $|M_{i+1}-M_i| \leq 1$ for every integer $i$. \\
If at step $i$, we have that $V(H)$ is not fully colored blue, then we must have $M_i = (\sum_{j=1}^i X_j) -(1- e^{-1})i$, and $\sum_{j=1}^i X_j\le \abs{S_i} < h$, as at each step $j$ for which $X_j=1$, at least one vertex was forced blue.
Therefore, 
\begin{align*}
    \bP(\Gamma' \ge 2h+t) &\leq \bP(M_{2h+t} \leq h - (1-e^{-1})(2h+t)) \leq \bP\left(M_{2h+t} \leq -\frac{t}{2}\right).
\end{align*}
By Azuma's inequality (see \Cref{Azuma inq}), and using $t \geq h$, we have 
\begin{align*}
    \bP\left(M_{2h+t} \leq -\frac{t}{2}\right) \leq \exp\left(- \frac{(t/2)^2}{2(2h+t)}\right)  \leq \exp\left(-\frac{t}{100}\right).
\end{align*}
% \mehdi{as $t\ge h$, we have $\exp\left(- \frac{(t/2)^2}{2(2h+t)}\right)  \leq \exp(\frac{-t^2}{24t})=\exp(\frac{-t}{24})$ as $2h+t\le 3t$. So we should change $\exp(\frac{-t}{100})$ with $\exp(\frac{-t}{24}).$}
This finishes the proof.
\end{proof}
%\rik{I'm trying to do the lemma below for arbitrary $s \in T_v$}
\begin{lemma}\label{T_v completely blue}
    Let $G$ be a connected graph on $n\ge 2$ vertices, and let $v \in V(G)$. 
    Let $(T_v,C)$ be a $v$-good pair, and let $s \in V(T_v)$. 
    Then for $n$ sufficiently large and any $t \geq 8\sqrt{n}$, after $8\sqrt{n}+t$ steps of probabilistic zero forcing on $G$ with initial blue set $S = \{ s \}$,  we have that $T_v$ is fully colored in blue with probability at least $1-\beta^t$, where $0 < \beta < 1$ is some absolute constant. 
   % \mehdi{On utilise hypo core seulement pour apply Lemma 4.5 no?}
%TO DO 25 05: utilises explicitement le mot modif pzf ds la pf.

    %\mehdi{need $s$ to be in core. in pf of th at the end, need vertices in $S$ also in core. mais c le cas car les $T_w$ sont rooted at w}
\end{lemma}
\begin{proof}
%We first remark that by \Cref{lemma 2.11 for Sun}, it suffices to prove this bound for the alternative probabilistic process where only vertices in $C$ are allowed to propagate their color (still independently and with same probabilities as in the probabilistic zero forcing).
Let $\Gamma=pt_{pzf}(G,\{s\},V(T_v))$ denote the random variable associated with the number of steps taken by the probabilistic zero forcing process on $G$, starting with the initial blue set $\{s\}$, until all the vertices in $T_v$ are colored blue.
%Let $\Gamma$ be the random variable associated to the time this alternative propagation lasts until all of $T_v$ is blue.
Since $(T_v,C)$ is $v$-good, there exists some $w\in C$ and $W\subseteq N(w)$ such that $\abs{V(T_v)\setminus (W \cup \{w\}}<|V(T_v) \setminus W|< 2\sqrt{n}$ and $C \setminus \{w\}$ is a core of $G[V(T_v) \setminus (W \cup \{w\})]$ in $G \setminus \{w\}$. 
% \mehdi{where do you use in the proof that $C\setminus \{w\}$ is a core of $G[V(T_v) \setminus (W \cup \{w\})]$ in $G \setminus \{w\}$?}
% \julien{It is used just after the equation \eqref{eq:Bound-Phase2}}
We first define the following stopping times: %\julien{Makes this more formal: write this with $\tau$}
\begin{enumerate}
    \item Let $\tau_1$ denote the first time at which some vertex in $C$ turns blue.
    \item Let $\tau_2$ be the first time after $\tau_1$ (i.e. $\tau_2 \geq \tau_1$) at which at least one vertex in $W \cup \{w\}$ is colored blue.
    % \mehdi{While I agree that $\tau_3\ge \tau_2$ since $N[w]$ being fully blue implies $W\cup \{w\}$ is blue, I think it's not correct that $\tau_2\ge \tau_1$. Indeed, we may (starting with a blue $s\in V(T_v)\setminus C$) first color a vertex in $W\cup \{w\}$ before "hitting" $C$. In this case, $\tau_2<\tau_1$.}
    % \julien{We define $\tau_2$ to be such that $\tau_2 \geq \tau_1$, so no problem on that.}
    \item Let $\tau_3$ be the first time after $\tau_2$ (i.e. $\tau_3 \geq \tau_2$) at which all vertices in $N[w]$ have been forced blue.
    %\item Let $\tau_3$ be the first time $\tau_2$ color entirely in blue the connected subgraph $T_v\setminus G[W]$ (if not done in (1)) by \Cref{lem:concentration of each T_v}, via core probabilistic zero forcing.
\end{enumerate}
Let $Z_1=\tau_1$, $Z_2=\tau_2-\tau_1$, $Z_3=\tau_3-\tau_2$ and $Z_4 = \Gamma - \tau_3$.
Intuitively, we are splitting the propagation into four successive phases: Phase 1 lasts until at least one vertex in $C$ turns blue, Phase 2 lasts until at least one vertex in $W \cup \{w\}$ is colored blue, Phase 3 ends when all of $N[w]$ is colored blue, and Phase 4 ends when all of $V(T_v)$ is blue, and $Z_1$, $Z_2$, $Z_3$ and $Z_4$ are respectively the random variables associated with the number of time steps that Phases 1, 2, 3 and 4 last.
%Let $T_1$ be the random variable which denotes the first time some vertex in $W$ becomes blue while running probabilistic zero forcing with initial blue set $\{s\}$ in $T_v$.
Suppose $t \geq 2 \sqrt n$.
If $s \in C$, then $\mathbb{P}(Z_1 \geq t)=0$. On the other hand, if $s \in V(T_v) \setminus C$, then $s$ has a neighbor in $C$. 
By \Cref{lemma 2.11 for Sun} we thus have that $Z_1$ is stochastically dominated by $pt_{pzf}(S_{d(s)},\{w_0\})$, where $S_{d(s)}$ is a star on $d(s)+1$ vertices centered at a vertex $w_0$ (we recall that for a graph $H$ and a vertex $u \in V(H)$, we defined $pt_{pzf}(H,\{u\})$ as the random variable associated with the time taken by the probabilistic zero forcing process on $H$, with only the vertex $u$ initially blue, until $H$ is fully blue).
% By \Cref{lemma 2.11 for Sun} we thus have that $\mathbb{P}(Z_1 \geq t )$ is at most $\mathbb{P}(pt_{pzf}(S_{d(s)},\{w_0\}) \geq t)$, where $S_{d(s)}$ is a star on $d(s)+1$ vertices centered at a vertex $w_0$ (we recall that for a graph $H$ and a vertex $u \in V(H)$, we defined $pt_{pzf}(H,\{u\})$ as the random variable associated with the time taken by the probabilistic zero forcing process on $H$, with only the vertex $u$ initially blue, until $H$ is fully blue). 
We have $t \geq 2 \sqrt n =\omega(\log n)$, so for $n$ sufficiently large, $\mathbb{P}(pt_{pzf}(S_{d(s)},\{w_0\}) \geq t) \leq \alpha^t$, for some $0 < \alpha < 1$, by Lemma \ref{lem: coloring d neighbours in a star}. Therefore,
\begin{align}\label{bound Z1}
    \mathbb{P}(Z_1 \geq t) \leq \alpha^t.
\end{align} 
% \mehdi{It is not clear why “Then $\mathbb{P}(Z_1 \ge t)$ is at most $\mathbb{P}(pt_{pzf}(S_0,{w_0}) \ge t)$.” Shouldn’t it be $\mathbb{P}(\text{$w_0$ becomes blue after $t$ steps}) \le \mathbb{P}(pt_{pzf}(S_0,w_0) \ge t)$?}
% \julien{It is indeed, but since $\mathbb{P}(Z_1 \geq t ) \leq \mathbb{P}(\text{$w_0$ becomes blue after $t$ steps})$ we are done.}
We now prove a bound for the tail of $Z_2$. Since some vertex $u$ in the core is blue at the end of Phase 1, as before, we consider the alternative probabilistic process where only vertices in $C$ propagate their color (independently and with the same probabilities as in the probabilistic zero forcing process on $G$). Let $Z_{2}'$ denote the number of steps taken by this alternative process, starting with only $u$ initially blue, until at least one vertex in $W \cup \{w\}$ is colored blue. Then we have $\mathbb{P}(Z_2 \geq 4 \sqrt n +t) \leq \mathbb{P}(Z_2' \geq 4 \sqrt n +t)$. 
We proceed similarly as in the proof of \Cref{lem:concentration of each T_v}.
%\julien{This is slightly annoying that we essentially have to repeat this proof, but I am worried this will make it more confusing for the reader if we merge the two in a more general and technical statement.}
Let $S_i$ be the set of blue vertices in $T_v$ at the $i$th step of the alternative probabilistic process.
Let $(X_i)_{i=1}^{\infty}$ be a sequence of random variables defined in the following way. Let $v^i$ be the vertex of minimum index in $C$ such that $v^i \in S_i$ and $N(v^i) \not\subseteq S_i$
(in other words, at step $i$, $v^i$ is blue and has at least one white neighbor in $T_v$).
Then we set $X_{i+1}=1$ if $v^i$ propagates its color to at least one of its white neighbors at step $i+1$, and $X_{i+1}=0$ otherwise.
If there exists no such $v^i$, we also set $X_{i+1}=0$. 
As in the proof of \Cref{lem:concentration of each T_v}, we remark that if at step $i$, no such $v^i$ as above exists, then all of $T_v$ is colored blue, and that if $S_i \neq V(T_v)$, we have
\begin{align*}
    \bP(X_{i+1} = 1 | S_0, \dots, S_i) \geq 1-e^{-1}.
\end{align*}
We define a sequence of random variables $(M_i)_{i=0}^{\infty}$ as follows. Set $M_0:=0$,  $M_{i+1}=M_i+X_{i+1}-(1-e^{-1})$ if there is some white vertex in $T_v$ at step $i$ and $M_{i+1}=M_i$ otherwise. As before, it follows that $(M_i)_{i = 0}^{\infty}$ is a submartingale with respect to $(S_i)_{i= 0}^{\infty}$, and that almost surely $|M_{i+1}-M_i| \leq 1$ for every integer $i$.
If at step $i$, we have that no vertex in $W \cup \{w\}$ is colored blue, then we must have $M_i = (\sum_{j=1}^i X_j) -(1- e^{-1})i$, and $\sum_{j=1}^i X_j \leq |V(T_v) \setminus (W \cup \{w\})| < 2\sqrt{n}$.
Therefore, 
\begin{align*}
    \bP(Z_2' \geq 4\sqrt{n}+t) \leq \bP( M_{4\sqrt{n}+t} \leq 2\sqrt{n} - (1- e^{-1})(4\sqrt{n}+t) ) \leq \bP\left( M_{4\sqrt{n}+t} \leq -\frac{t}{2} \right).
\end{align*} %\mehdi{first inequality not clear}
By Azuma's inequality (see \Cref{Azuma inq}), and assuming $t \geq 2 \sqrt{n}$, we have 
\begin{align*}
    \bP\left(M_{4\sqrt{n}+t} \leq -\frac{t}{2}\right) \leq \exp\left(- \frac{(t/2)^2}{2(4\sqrt{n}+t)}\right)  \leq \exp\left(-\frac{t}{100}\right).
\end{align*}
We can then conclude that 
\begin{align}
\label{eq:Bound-Phase1}
    \bP(Z_2 \geq 4\sqrt{n}+t) \leq \exp\left(-\frac{t}{100}\right).
\end{align}
We now turn to $Z_3$.
We condition on the vertex of smallest index $w^h$ that first became blue in $W \cup \{w\}$.
We further split Phase $3$ into two subphases: let $\tau$ be the first time that the vertex $w$ becomes blue, and let $Z^a_3=\tau-\tau_2$ and $Z^b_3=\tau_3-\tau$.
It follows that $Z_3=\tau_3-\tau_2=Z^a_3+Z^b_3$, and thus %\julien{Add conditioning?}
\begin{align*}
    \bP( Z_3 \geq t | w^h) \leq \bP(Z^a_3 \geq t/2 | w^h)+\bP(Z^b_3 \geq t/2|w^h).
\end{align*}
Conditioned on $w^h$, the random variables $Z^a_3$ and $Z^b_3$ are respectively stochastically dominated by $pt_{pzf}(S_a,\{w_a\})$ and $pt_{pzf}(S_b,\{w_b\})$, where $S_{a}$ is a star on $d(w^h)+1$ vertices centered on a vertex $w_a$ and $S_{b}$ is a star on $d(w)+1$ vertices centered on a vertex $w_b$.
%\julien{A bit more explanation for this is probably needed: maybe conditioning}
Thus by \Cref{lem: coloring d neighbours in a star}, and assuming $t \geq 2\sqrt{n} = \omega(\log n)$, for $n$ sufficiently large, we have
\begin{align*}
   \bP( Z_3 \geq t | w^h) \leq \bP(Z^a_3 \geq t/2 | w^h) + \bP(Z^b_3 \geq t/2 | w^h) \leq 2\alpha^{t/2}.
\end{align*}
% \mehdi{for the last inequality say we use Lemma 3.8 before applying Lemma 3.9}
% \julien{It is true that we are using Lemma 3.8 many times, and we do not say everywhere. But I think at this point of the paper, it is ok if we omit it.}
As this upper bound does not depend on $w^h$, we obtain by summing over all possible outcomes of $w^h$ that
\begin{align}
\label{eq:Bound-Phase2}
   \bP( Z_3 \geq t) \leq 2\alpha^{t/2}.
\end{align}
Finally, as $(T_v,C)$ is $v$-good, we have that $C \setminus \{w\}$ is a core of $G[V(T_v) \setminus (W \cup \{w\})]$ in $G \setminus \{w\}$ and that $|V(T_v) \setminus (W \cup \{w\})| < |V(T_v) \setminus W|< 2\sqrt{n}$. As $C$ is connected, and $w \in C$, at least one vertex in $C \setminus \{w\} $ is blue at the end of Phase 3 (since $N(w)$ is colored blue). As $w$ is already blue at the end of Phase 3, the probability that a vertex $w' \in V(G) \setminus \{w\}$ forces its white neighbors blue in the probabilistic zero forcing process on $G \setminus \{w\}$ is at most the corresponding probability in $G$.

Therefore, by \Cref{lemma 2.11 for Sun,lem:concentration of each T_v} for the graph $G[V(T_v) \setminus (W \cup \{w\})]$ with core $C \setminus \{w\}$, and assuming $t \geq 2\sqrt{n}$, we obtain
\begin{align}
\label{eq:Bound-Phase3}
    \bP(Z_4\ge 4\sqrt{n}+ t)&\le \exp\left(-\frac{t}{100}\right). 
\end{align} 
% \mehdi{need to say WMA $u\in V(T_v)\setminus(W\cup \{w\})$, otherwise Phase 2 not needed. This assumption on $u$ is important.}
% \mehdi{Need to say WMA $V(T_v)\setminus (W\cup \{w\})$ not fully blue at the end of Phase 3}
% \julien{It is not compulsory to say it, as otherwise the relevant Phase lasts for 0 step, so not an issue}
By putting together \eqref{bound Z1}, \eqref{eq:Bound-Phase1}, \eqref{eq:Bound-Phase2} and \eqref{eq:Bound-Phase3}, we find that 
\begin{align*}
    \bP( Z_1+Z_2+Z_3&+Z_4 \geq 8\sqrt{n}+ 4t) \\
    &\leq \bP(Z_1 \geq t)+\bP(Z_2\ge 4\sqrt{n}+ t)+\bP(Z_3\ge  t)+\bP(Z_4\ge 4\sqrt{n}+ t) \\
    &\leq \alpha^t + \exp\left(-\frac{t}{100}\right)+2\alpha^{t/2}+\exp\left(-\frac{t}{100}\right),
\end{align*}
% as the probability of the event $Z_4\ge \frac{t}{2}$ for the probabilistic zero forcing is at most the probability of this same event for the modified process given in \Cref{lemma 2.11 for Sun}.
and thus, there exists some absolute constant $\beta <1$ such that
\begin{align*}
    \bP(\Gamma \geq  8\sqrt{n}+ 4t) \leq \beta^{4t}.
\end{align*}
By rescaling $t$, we obtain the desired conclusion.
% On the other hand, we have by running core probabilistic zero forcing, %\julien{This application of the lemma is not correct: you cannot just ignore $W$ like this. You need to justify why you can apply this lemma. Essentially it is because we can "remove the blue vertices in $W$" once $W$ is fully colored, but this needs to be made rigorous.}
% \begin{align*}
% \bP(Z_3 \ge 4\sqrt{n} + t/2) 
% &\le \bP(Z_3 \ge 2\abs{V(T_v) \setminus W} + t/2) \le  10 \exp\left(-\frac{t}{200}\right)
% \end{align*} where \Cref{lem:concentration of each T_v} applied to $T_v\setminus G[W]$ (we may assume $t>4\sqrt{n}$ in order to apply this lemma, but for simplicity we keep $t>2\sqrt{n}$) justifies the second inequality since \Cref{lemma 2.11 for Sun} implies that once $W$ is entirely colored in blue, the probability that $Z_3\ge \abs{V(T_v)\setminus W}$ lower bounds the probability of this event but if we followed the core probabilistic zero forcing. 
% As $Z_2\le Z_3+Z_4$, for $t>2\sqrt{n}$, we conclude that  
% \begin{align*}
%     \bP(Z_2\ge 8\sqrt{n}+t)\le \bP(Z_3\ge 4\sqrt{n}+\frac{t}{2})+\bP(Z_4\ge \frac{t}{2})\le  10 \exp\left(-\frac{t}{200}\right)+2\alpha^{t/4}\le  12 \exp\left(-\frac{t}{100}\right)
% \end{align*}
%  where we are able to apply the inequality derived for $\bP(Z_3\ge 4\sqrt{n}+t/2)$ since the probability of this event for the normal process is bounded above by the probability for the core process.
% Then we conclude that after $8\sqrt{n}+t$ steps, for $t>2\sqrt{n}$, we obtain that $T_v$ is completely blue with probability at least $1- 12 \exp\left(-\frac{t}{100}\right)$
\end{proof}

We are now ready to prove \Cref{thm:Throttling}.
\begin{proof}[Proof of Theorem~\ref{thm:Throttling}.]
Let $G$ be a connected graph on $n$ vertices, and note that we may assume $n \geq 4$. 
Label the vertices in $V(G)$ by $v_1,\dots,v_n$. 
For every $i$, we let $(T_{v_i},C_{v_i})$ be a $v_i$-good pair obtained by applying Lemma \ref{decomposition core leaves}. 
Let $S$ be a subset of $V(G)$ of size $\lfloor 2\sqrt{n}\rfloor$ given by \Cref{lem:Set-System} applied to the sets $V(T_{v_1}), \dots, V(T_{v_n})$, and denote by $\Gamma$ the random variable associated with the time taken by the probabilistic zero forcing process on $G$, starting with the set $S$ initially blue, until $G$ is fully colored blue.
Then it suffices to prove that $\mathbb{E}[\Gamma]=O(\sqrt{n})$. 
We let $\tau_1$ be the first time that all the vertices in $\cup_{v \in S}T_{v}$ are colored blue.
We then let $Z_1=\tau_1$ and $Z_2=\Gamma-\tau_1$.
Applying \Cref{T_v completely blue} for each $(T_v,C_v)$ with $v \in S$ and performing a union bound over all $v \in S$, we find that for all $t \geq 8 \sqrt{n}$, we have
\begin{align}
\label{eq:Bound-Phase1-MainProof}
\bP(Z_1 \geq 8\sqrt{n}+t) = O(\sqrt{n}\beta^t).
\end{align}
Note that at time $\tau_1$, each $T_v$ for $v \in V(G)$ has at least one blue vertex.
Thus, once again applying \Cref{T_v completely blue} for each $(T_v,C_v)$ with $v \in V(G)\setminus S$ and performing a union bound over all $v \in V(G)\setminus S$, we find that for all $t \geq 8\sqrt{n}$, 
\begin{align}
\label{eq:Bound-Phase2-MainProof}
\bP(Z_2 \geq 8\sqrt{n}+t) = O(n\beta^t).
\end{align}
Combining \eqref{eq:Bound-Phase1-MainProof} and \eqref{eq:Bound-Phase2-MainProof}, we get for all $t \geq 8\sqrt{n}$, 
\begin{align*}
    \bP(\Gamma \geq 16\sqrt{n}+2t) &\leq \bP(Z_1 \geq 8\sqrt{n}+t)+ \bP(Z_2 \geq 8\sqrt{n}+t) =O(n \beta^t).
\end{align*}
%It follows that
%\begin{align*}
    %\mathbb{E}[\Gamma1_{\Gamma > 28\sqrt{n}}] \leq \sum_{t \geq \sqrt{n}} \bP(\Gamma \geq 16\sqrt{n}+2t) (16\sqrt{n}+2t+2) =O(\sum_{t \geq \sqrt{n}} t^3\beta^t) = o(1),
%\end{align*}
%and thus
%\begin{align*}
    %\mathbb{E}[\Gamma] &= \mathbb{E}[\Gamma1_{\Gamma \leq 28\sqrt{n}}] + \mathbb{E}[\Gamma1_{\Gamma > 28\sqrt{n}}] = O(\sqrt{n}),
%\end{align*}
%as wanted.
Thus, we have
\begin{align*}
    \sum_{k > 32 \sqrt n} \mathbb{P}(\Gamma \geq k)= O\left(\sum_{t \geq 8 \sqrt n}n 
    \beta^t\right)=O\left(\sum_{t \geq 8 \sqrt n}t^2 
    \beta^t\right)=O(1).
\end{align*}
From the tail sum formula for expectation, it follows that 
\begin{align*}
    \mathbb{E}[\Gamma]=\sum_{k=1}^{\infty} \mathbb{P}(\Gamma \geq k) &= \sum_{k \leq 32 \sqrt n} \mathbb{P}(\Gamma \geq k) + \sum_{k > 32 \sqrt n} \mathbb{P}(\Gamma \geq k)= O(\sqrt n)+O(1)=O(\sqrt n), 
\end{align*}
as desired.

\end{proof}
\section{Upper bound on the expected propagation time of probabilistic zero forcing} \label{section ept improved}

%\julien{Rik: here you write the lemma "we have $\mathbb{E}[\min(X_{i+1}, X_i+100)] \geq X_i+2$ unless... "}
%As in \cite{Sun}, the following structure will play a crucial role in our analysis.
A pair \( (v, w) \) of vertices in a graph \( G \) is said to be a \emph{cornerstone} if $\{v,w\}$ is a cut-set of $G$, and \( v \) and \( w \) share a common neighbor in \( G \). 
% \mehdi{or are connected. see definition of Narayanan and Sun} \julien{I think that Rik managed to rule out this case in his proof of \ref{lem: two new blue vertices on average, after truncation}}
Given a cut-set $A \subseteq V(G)$, a pair of disjoint subsets $(S,T)$ is called a \emph{valid pair} if $S \cup T = V(G) \setminus A$, there are no edges between $S$ and $T$, and $|S| \leq |T|$.

As explained in \Cref{sec:Outline}, we now define a function for each cut-vertex and cornerstone, according to which we will make the choice of the initial blue vertex.
%and we will take as starting point a cut-vertex or cornerstone which maximises this carefully chosen function. 
% This formalises the intuition of selecting a starting vertex "near the centre of the graph". \julien{Find better transition}

%by proposing a new selection of the initial cut-vertex or cornerstone that optimises a specifically designed function.
\begin{definition}\label{new_choice_of_S_T}
    If $v$ is a cut-vertex, define $h(v)$ to be the maximum of $|S|-d_S$ over all valid pairs $(S,T)$ with respect to the cut-set $\{ v\}$, where $d_S := |S \cap N[v]|$. 
    % \mehdi{it should be $d_S := |S \cap N(v)|$ as $v\notin S$.} \julien{It is better how it is now because it is then uniform with the definition for the cornerstone} %If \( v \) is not a cut-vertex, set $h(v)=0$. Note that $h(v)\ge 0$ for any vertex $v$, since $d_S \leq |S|$. 
    
    Likewise, if \( (v, v') \) is a cornerstone, we define $h(v,v')$ to be the maximum of $|S| -d_S$ over every common neighbor $y \in N(v) \cap N(v')$ and valid pairs $(S,T)$ with respect to the cut-set $\{ v,v'\}$, where $d_S := |S \cap N[\{v,v',y\}]|$. %\mehdi{how do you choose $y$ if $(v,v')$ share many common neighbors? should we maximise over $y$ such that $h(v,v')$ is maximal?} 
\end{definition}

    Let $G$ be a graph. In a probabilistic zero forcing process on $G$, we say that one \emph{encounters} a cut-vertex $v$ or a cornerstone $(v,v')$ at step $i$ if one of the following occurs at step $i$: %\julien{Try to make this definition shorter}
    \begin{itemize}
        \item $v$ is the only blue vertex that is adjacent to some white vertex. 
        \item $v$ is the only white vertex that is adjacent to some blue vertex. 
        \item $v,v'$ share a common white neighbor and are the only two blue vertices that are adjacent to some white vertex.
        \item $v,v'$ share a common blue neighbor and are the only two white vertices that are adjacent to some blue vertex.
        %\item  
        %there are at least three blue vertices that are adjacent to some white vertex in $G$, and that exactly two white vertices are adjacent to any blue vertex in $G$. Furthermore, $v,v'$ are white and there exists a blue vertex adjacent to them.
        %\item 
%there are at least three blue vertices denoted $v',b_2,\dots,b_k$ (with $k\ge 3$) adjacent to white vertices denoted $v,w_2,\dots,w_k$ and at least three white vertices adjacent to blue vertices. 
%All blue vertices except $v'$ are adjacent to exactly one white vertex, and $v$ is a white vertex adjacent to all blue vertices. 
%In this configuration, $v'$ is adjacent to all white vertices $v, w_2, \ldots, w_k$ that have any blue neighbors, and none of $b_2, \ldots, b_k$ are adjacent to any blue vertex other than $v'$. \rik{don't think I need this one for the next lemma} \julien{ok great}
    \end{itemize}

%We introduce the following notation: %\julien{Remove those notations and simply write $d(v)$, $d_S(v)$ and $d_T(v)$}
%$d:=\deg(v)$, $d':=\deg_S(v)$, $d'':=\deg_T(v)$ and
% \[
% \alpha_S := d',\quad \max(d', \deg_S(v')),\quad \text{or} \quad \max(d', \deg_S(v'), \deg_S(v''))
% \]

% \[
% \alpha_T := d'',\quad \max(d'', \deg_S(v')),\quad \text{or} \quad \max(d'', \deg_S(v'), \deg_S(v''))
% \]

% \[
% \log(\alpha_S) := \log d',\quad \max(\log d', \log \deg_S(v')),\quad \text{or} \quad \max(\log d', \log \deg_S(v'), \log \deg_S(v''))
% \]

% \[
% \log(\alpha) := \log d \quad \max(\log d, \log \deg(v')),\quad \text{or} \quad  \max(\log d,\log \deg(v'),\log \deg(v''))
% \]
% depending on whether a single vertex $v$ was selected,
% or a pair $v, v'$ that share an edge, or a pair that share a common neighbor $v''$.

We first prove the following lemma, which is along similar lines as Lemma 4.4 in \cite{Sun}.

\begin{lemma}
    \label{lem: two new blue vertices on average, after truncation} Let $X_{i}$ be the random variable associated with the number of blue vertices in $G$ after $i$ steps of probabilistic zero forcing. Suppose that every blue vertex in $G$ has at least one blue neighbor after the $i$th step of probabilistic zero forcing. Then, unless $G$ is fully blue or we encounter a cut-vertex or cornerstone at the $i$-th step, we have
    \begin{equation*}
        \mathbb{E}[\min(X_{i+1}, X_{i}+25)] \geq X_{i} + 2.
    \end{equation*}
\end{lemma}

\begin{proof}
    Let $Y_{i+1} := X_{i+1}-X_{i}$. We want to show that $\mathbb{E}[\min(Y_{i+1},25)] \geq 2.$ Note that $Y_{i+1}$, the number of vertices which turn blue in the $(i+1)$-th step, can be viewed as the sum of independent Bernoulli random variables. 
    Therefore, it suffices to show that either $\mathbb{E}[Y_{i+1}] \geq 2+1/5$, in which case we would be done by Lemma \ref{lem: high expectation after truncation}, or that $\mathbb{E}[\min(Y_{i+1},25)] \geq 2.$  %\julien{Say at the beginning that if $\mathbb{E}[Y_{i+1}] \geq 2+1/5$ then by Lemma \ref{lem: high expectation after truncation}, we are done, so we focus on proving either. So that we do not have to use this lemma repeatedly.} 
    Note that a blue vertex with $w$ white neighbors will force each of its white neighbors blue independently with probability at least $2/(w+1)$, since every blue vertex has at least one blue neighbor. Let us denote by $v_{1}, \dots , v_{k}$ the blue vertices in $G$ with at least one white neighbor, and by $w_{1}, \dots , w_{l}$ the white vertices in $G$ with at least one blue neighbor after $i$ steps. Since we do not encounter a cut-vertex, we may assume that $k,l \geq 2$. If $k=2$, then $v_{1}$ and $v_{2}$ do not have a common white neighbor, as otherwise this would contradict the assumption that we do not encounter a cornerstone. If $l=2$, then both $v_1$ and $v_2$ have precisely one white neighbor, so $\mathbb{E}[Y_{i+1}]=\mathbb{E}[\min(Y_{i+1},25)]=2$. If $l \geq 3$, then either $v_{1}$ or $v_2$ has at least two white neighbors. Assume without loss of generality that $v_1$ has $l_1 \geq 2$ white neighbors. Then $v_1$ will force at least $2l_1/(l_1+1)\geq1+1/3$ of its white neighbors blue in expectation, so $\mathbb{E}[Y_{i+1}]\geq 2 + 1/3$. %\julien{Add a small computation for this}, and we are done by Lemma \ref{lem: high expectation after truncation}. 
    If $k \geq 3$ and $l=2$, then $w_1$ and $w_2$ do not share a blue neighbor, as otherwise this would contradict the assumption that we do not encounter a cornerstone. 
    Then each blue vertex is adjacent to either $w_1$ or $w_2$, and so both $w_1$ and $w_2$ will be forced blue with probability 1. Hence, $\mathbb{E}[Y_{i+1}]=\mathbb{E}[\min(Y_{i+1},25)]=2$. We may henceforth suppose that $k,l \geq 3$. Let $a:=\text{deg}_w(v_1), b:=\text{deg}_w(v_2)$ and $c:=\text{deg}_w(v_3)$, where $\text{deg}_w(v_{i})$ denotes the number of white neighbors of $v_i$, and assume without loss of generality that $a \geq b \geq c$. If $a \geq 3$ and $b \geq 2$ %\julien{This follows because the formula below is increasing in $a$ and $b$ and $c$ right? If yes, then say it}
    , we have 
    \begin{equation*}
        \mathbb{E}[Y_{i+1}]\geq a\cdot\frac{2}{a+1}+b\cdot\frac{2}{b+1}\left(1-\frac{2}{a+1}\right)+c\cdot\frac{2}{c+1}\left(1-\frac{2}{a+1}\right)\left(1-\frac{2}{b+1}\right) \geq 2 + \frac{1}{3},
    \end{equation*}
    where the last inequality holds since the expression above is increasing in $a,b$ and $c$. If $a \leq 2$, then each $v_j$ has at most two white neighbors. Thus, each $w_j$ will be forced blue with probability at least $2/3$. 
    Therefore $\mathbb{E}[Y_{i+1}] \geq l\cdot \frac{2}{3}$, and if $l \geq 4$, we then have $\mathbb{E}[Y_{i+1}] \geq 2+\frac{2}{3}$.
    On the other hand, if $l=3$, then $\mathbb{E}[Y_{i+1}]=\mathbb{E}[\min(Y_{i+1},25)]\geq 2$. 
    Finally, suppose that $b=1$. Then one of the white vertices will be forced blue with probability 1, while the other $l-1$ white vertices will be forced blue with probability at least $2/(l+1)$. Thus, $\mathbb{E}[Y_{i+1}]\geq 2 \left(\frac{l-1}{l+1} \right)+1 = 3-\frac{4}{l+1}$. If $l \geq 4$, then $\mathbb{E}[Y_{i+1}] \geq 2+1/5$, and if $l=3$, then $\mathbb{E}[Y_{i+1}]=\mathbb{E}[\min(Y_{i+1},25)]\geq 2$. This proves the lemma.
\end{proof}

In the spirit of the proof of Narayanan and Sun \cite{Sun}, we provide the following algorithm.
\begin{enumerate}\label{new algo for ept}
    \item Choose a single vertex \( v \) or a cornerstone \( v, v' \) such that the value of \( h(v) \) or \( h(v, v') \) in Definition \ref{new_choice_of_S_T} is maximized. \\
    In other words, if $\max_{w \text{ cut-vertex}}h(w)\geq \max_{(w,w') \text{ cornerstone}}h(w,w')$, then we choose $v$ such that $h(v)=\max_{w \text{ cut-vertex}}h(w)$. Otherwise, we choose $v,v'$ such that $h(v,v')=\max_{(w,w') \text{ cornerstone}}h(w,w')$.
    \\
If we chose a cut-vertex $v$,
pick a valid pair \( (S,T) \) which maximises $h(v)$.
%\julien{It is important that the sets $(S,T)$ chosen are the ones maximising the definition.} such that \( S \cup T = V \setminus \{v\} \), there are no edges between \( S \) and \( T \) and $\left|S\right|\leq \left|T\right|$.
    Likewise, if we chose a cornerstone, pick a valid pair $(S,T)$ (as well as a common neighbor $y$ of $v$ and $v'$) which maximises $h(v,v')$.  
    If there are no cut-vertices or cornerstones, choose any $v\in V$, and take \( S = \emptyset \) and $T= V(G) \setminus \{v\}$.  %\julien{Is that important?}
    Let $d_S$ be the quantity defined in Definition \ref{new_choice_of_S_T} if we picked a cut-vertex or cornerstone, and let $d_S=0$ otherwise. Let $d_T=|T \cap N[v]|$ if we picked a single vertex $v$, and let $d_T= |T \cap N[\{v,v',y\}]|$ if we picked a cornerstone $v,v'$ with a common neighbor $y$.    
    
    \item Initialize with \( v \) blue and all other vertices white.
        \setcounter{enumi}{2}
    \item Pick some arbitrary ordering \( v_1, \dots, v_n \) of the vertices, where \( v, v' \) can be labeled with any number.
    \item Run the probabilistic zero forcing process until all neighbors of \( v \) are blue, and if we picked a cornerstone \( v, v' \) in Step 1, until all neighbors of \( v' \) and $y$ are also blue. 
    %\item Turn white all vertices that are not \( v, v' \), or any of their neighbors. 
\item Run the probabilistic zero forcing process on the induced subgraph \( G[T] \), until the number of white vertices in \( T \) is at most
$\abs{S}-d_S$. 
% \mehdi{shouldn't it be at most
% $\abs{S}-d_S+3$ for the same technicality reasons that N. and S. considered $\abs{S}+3$ in their definition of Step 5? (the $+3$ should count for the case we have a cornerstone i.e. the vertices $v,v',y$.)}
% \julien{To me the proof is correct as it is}
%\item If we have encountered a cut-vertex or a cornerstone in the previous step, then we first continue the propagation until at most $|S|-d_S$ white vertices are left.
%\julien{Say that if we encounter a cut-vertex or cornerstone, then we continue until the cut-vertex or cornerstone is blue, and then we stop. Actually we should put this in Step 6.}

%\julien{Add a new step 6}

    \item At each step, suppose there is some blue vertex in \( G[T] \) with $k\ge1$ white neighbors in \( G[T] \). Choose such a blue vertex \( v_i \in G[T] \) with smallest index \( i \), and run probabilistic zero forcing but where each white neighbor of \( v_i \) becomes blue with probability 1 if \( v_i \) only has exactly one white neighbor, and becomes blue with probability \( \frac{4}{3k} \) otherwise.
    If no such blue vertex exists in $G[T]$, do nothing.
    Apply the same procedure on $G[S]$ in parallel.
    %\julien{I think we need to add that if we encountered something, then we start by doing a star propagation.}
\end{enumerate}

As mentioned in the outline of our proofs, this algorithm is similar to the one in Narayanan and Sun's proof \cite{Sun}, where the key difference arises from our choice of the function $h$ to maximise in Step 1, which also modifies the analysis needed to bound the expected runtime of Step 5.

We remark that Step 6 does not decrease the total expected propagation time:
indeed, letting $b$ and $k$ be, respectively, the number of blue and white vertices in the neighborhood of $v_i$, we note that for $k \geq 2$, we have
\begin{align*}
   \frac{b+1}{b+k} \geq \frac{4}{3k},
\end{align*}
and thus by \Cref{lemma 2.11 for Sun}, Step 6 does not decrease the total expected propagation time.
 Therefore, it suffices to bound the expected runtimes of Steps 4, 5, and 6 and show that their sum is at most $\frac{n}{2}+O(1)$, which will imply the desired bound on the expected propagation time. 
 %Each of them can be obtained via modifications to Lemmas 4.8, 4.9 and 4.10 in \cite{Sun}.
 %Therefore, we state those bounds here and defer their proof to the Appendix. 

 \begin{lemma}
    \label{lem:Bounds-Step4}
    Step 4 takes $O(\log (d_T+d_S))+O(1)$ steps in expectation.
\end{lemma}

The proof of Lemma \ref{lem:Bounds-Step4} is similar to that of \cite[Lemma 4.8]{Sun}, and can be found in \Cref{Appendix-Bounds-Step4}.

%\julien{Write the proof of the next lemma here: there is a distinction of cases to do here: either we encounter a cutvertex or cornerstone, in which case we need to split in 3 subphases: until we encounter the cutvertex or cornerstone, then we colour to have the cutvertex or cornerstone blue, and then we do the star colouration to colour fast the number of vertices necessary to get to $|S|-d_S$ vertices. If we do not encounter cutvertex or cornerstone, then easier.}
 
\begin{lemma}
    \label{lem:Bounds-Step5}
    Step 5 takes at most $\frac{|T|-d_T-|S|+d_S}{2}+O(1)$ steps in expectation.
    %\julien{Add the expectation there}
\end{lemma}
\begin{proof}
In Step 5, we run the probabilistic zero forcing process on the induced subgraph $G[T]$ until we have at most $|S|-d_S$ white vertices left in $G[T]$. Moreover, during this step, each blue vertex has at least one blue neighbor. 
We can therefore use Lemma \ref{lem: two new blue vertices on average, after truncation}, along with Doob's Optional Stopping Theorem (see \Cref{thm:Doob}) to bound the expected running time of this step. Let $Z_5$ denote the random variable associated with the running time of Step 5.
Let $S_i$ denote the set of blue vertices in $G[T]$ after $i$ iterations of probabilistic zero forcing during Step 5, and let $X_i := |S_i|$. Note that $X_0 \geq d_T$, and let $\tau_1$ be the first time that there are at most $|S|-d_S$ white vertices left in $G[T]$, or we encounter a cut-vertex or cornerstone. Clearly, $\tau_1$ is a well-defined stopping time with respect to $(S_i)_{i=0}^{\infty}$. Finally, let $\tau_2 := Z_5 - \tau_1$. 
% \mehdi{WMA $G[T]$ not fully blue} \julien{Concerning this comment and all the above similar ones: I think we do not need to say this: if this was the case everything we say is still true, and Stage 5 lasts 0 steps, so we are all good.}
%\julien{This is not fully formal to condition on the future events (i.e. not encountering a cornerstone or cutvertex): maybe to be fully formal we can define $\tau$ as the stopping time when the above happens or we have encountered a cornerstone or cutvertex? And then merge the two first phases of the distinction of cases?}. 

We define $M_0:=X_0$, $M_{i}:= \min(X_{i}, X_{i-1}+25)-2i$ if we did not encounter a cut-vertex or cornerstone during the first $i-1$ iterations of probabilistic zero forcing during Step 5 and $M_{i}:=M_{i-1}$ otherwise. 
%More precisely, we let $X_0, X_1, \ldots$ be random variables where $X_i$ represents the number of blue vertices at time $i$. 
%Initially, we have $X_0 = \alpha_T$, and the corresponding set is denoted as $S_0$.

%Let $\tau$ be the random variable representing the number of steps until Step $5$ is done.
%In other words, $\tau$ is the smallest index $i$ for which $X_i = \abs{\widetilde{T}}$. 
%Clearly, $\tau$ is a well-defined stopping time. %\julien{This is bad notation: it feels like we are conditioning on the final outcome, which we should not. I will rewrite it.}

%Next, consider the sequence of random variables $M_0, M_1, \ldots$, defined by $M_n = X_n - 2n$. Observe that $M_0 = X_0 - 0 = \alpha_T$.
%As demonstrated above, conditionally on $i \leq \tau$, the expected number of blue vertices increases by $2$ over step $i+1$, that is
%\begin{equation}
    %\mathbb{E}[X_{i+1} \mid X_i] \geq X_i + 2.
    %\tag{1}\label{eq:one}
%\end{equation}
\noindent
We claim that $\{M_i\}_{i=0}^\infty$ forms a submartingale with respect to $\{S_i\}_{i=0}^\infty$. Indeed, if we did not encounter a cut-vertex or cornerstone during the first $i$ iterations, 
\begin{align*}
    \mathbb{E}[M_{i+1} \mid S_0, S_1, \ldots, S_i]
    &= \mathbb{E}[\min(X_{i+1},X_{i}+25)-2(i+1) \mid S_0, S_1, \ldots, S_i]\\ &\geq X_{i}+2-2(i+1) \\&\geq \min(X_i,X_{i-1}+25)-2i =M_i,
\end{align*}
where the first inequality follows from Lemma \ref{lem: two new blue vertices on average, after truncation}. Otherwise, we have $\mathbb{E}[M_{i+1} \mid S_0, S_1, \ldots, S_i]=M_i$. Additionally, the absolute difference $\abs{M_i - M_{i-1}}$ is almost surely uniformly bounded by $n+1$ for all $i \geq 1$. It is also easy to see that $\mathbb{E}[\tau_1] < \infty$.
Applying Doob’s Optional Stopping Theorem (see \Cref{thm:Doob}), we obtain
\[
\mathbb{E}[M_{\tau_1}] = \mathbb{E}[\min(X_{\tau_1},X_{\tau_1-1}+25)] - 2\mathbb{E}[\tau_1] \geq \mathbb{E}[M_0] \geq d_T.
\]
%Thus, we conclude that
%\[
%\text{ept}(G[\widetilde{T}], S_0) = \mathbb{E}[\tau] \leq \frac{{\mathbb{E}[\abs{T \setminus T'}]} - \alpha_T}{2}\le \frac{\abs{T}-\alpha_T-{\mathbb{E}[\abs{T'}]}}{2},\]
%as required.
% \mehdi{We used $M_{\tau_1}:= \min(X_{\tau_1}, X_{\tau_1-1}+25)-2\tau_1$ which holds if we did not encounter a cut-vertex or cornerstone during the first $\tau_1$ iterations of probabilistic zero forcing during Step 5. However, at step $\tau_1$ we may encounter a cut-vertex or a cornerstone. Thus, we may not have this equality at $\tau_1$ (but at $\tau_1-1$, it holds)}
% \julien{To me this is still correct, but the "i+1 iterations" above should be changed into "i iterations" }
Thus, we have 
\begin{align}
\label{eq: bound tau1}
\mathbb{E}[\tau_1]\leq \frac{\mathbb{E}[\min(X_{\tau_1},X_{\tau_1-1}+25)]}{2}-\frac{d_T}{2}.
\end{align}
Let $S'$ denote the set of white vertices in $G[T]$ at time $\tau_1$.  Let $A$ denote the event that there are at most $|S| - d_S$ white vertices in $G[T]$ at time $\tau_1$, and let $B$ denote the event there are more than $|S| - d_S$ white vertices in $G[T]$ at time $\tau_1$. 
Note that by definition of $\tau_1$, if the event $B$ occurs, then a cut-vertex or a cornerstone is encountered at time $\tau_1$.  \\
%From \Cref{eq: bound tau1}, it follows that 
%\begin{align*}
    %\mathbb{E}[\tau_1 \mid A] \leq  \frac{\mathbb{E}[X_{\tau_1-1}\mid A]}{2}-\frac{d_T}{2} +O(1) \leq \frac{|T|-d_T-|S|+d_S}{2}+O(1).
%\end{align*}
Observe that $\mathbb{E}[\tau_2 \mid A]=0$.
Next, we condition on the event $B$. Let us assume that in Step 1, we chose a cornerstone $(v,v')$ with common neighbor $y$. The other cases can be handled similarly. Let $B' \subseteq B$ denote the event that we encountered a cornerstone $(w,w')$ with common neighbor $z$ at time $\tau_1$, where $w,w'$ are the only blue vertices in $G[T]$ with at least one white neighbor and $z$ is white. 
%\mehdi{$z\in T$}
% \mehdi{what about the other case of encounter: $w,w'$ are white?} \julien{Yes they also need to be dealt with}
%\mehdi{maybe we should mention that if $z$ is blue then $\tau_2^a=0$, so we may assume $z$ to be white.}
%From \Cref{eq: bound tau1}, we get \julien{To me the next equation is not correct: we cannot simply condition on $B'$ from the previous equation. I think we can deduce it through the same proof, but I do not think this is needed, as you can do the conditioning on $\tau_2$ only.}
%\begin{align}
%\label{eq: bound tau1 encounter}
    %\mathbb{E}[\tau_1 \mid B'] \leq \frac{\mathbb{E}[X_{\tau_1-1} \mid B']}{2}-\frac{d_T}{2}+O(1) \leq \frac{|T|-|S_2|-d_T}{2}+O(1).
%\end{align}
%Now suppose that we encounter a cut-vertex or cornerstone during Step 5. We further assume that in Step 1, we chose a cornerstone $(v,v')$ with common neighbor $y$, and during Step 5, we encounter a cornerstone $(w,w')$ with common neighbor $z$, where $w,w'$ are the only blue vertices in $G[T]$ with at least one white neighbour, and $z$ is white. The other cases can be dealt with in a similar manner. In this case, we split Step 5 into three sub-phases. First, we run probabilistic zero forcing on $G[T]$ until we encounter the cornerstone $(w,w')$. Let $S_2$ denote the set of white vertices in $T$ at the end of the first sub-phase, and 
We define $T' := V(G)\setminus (S' \cup \{w,w'\})$, $d_{S'}:= |S' \cap N[\{w,w',z\}]|$ and $d_{T'}:= |T' \cap N[\{w,w',z\}]|$. %\mehdi{need to define $S'$} \julien{It is defined above}

\begin{claim}\label{claim: set inclusion in encounter}
    $S \setminus N[\{v,v',y\}] \subsetneq T' \setminus N[\{w,w',z\}]$. Consequently, $|T'|>|S'|$.
\end{claim}

\begin{proof}
 Since $G$ is connected, $N(\{v,v'\}) \cap S \neq \emptyset$. Thus, it suffices to show that $S \subseteq T' \setminus N[\{w,w',z\}]$. We have $S \cap (S' \cup \{w,w'\}) = \emptyset$, which implies that $S \subseteq T'$. Note that $w,w',z \notin S \cup \{v,v'\}$, as $v,v'$ do not have any white neighbors after Step 4. Therefore, $N[\{w,w',z\}] \cap S = \emptyset$, so we have $S \subseteq T' \setminus N[\{w,w',z\}]$. %\julien{This proof is very concise actually, nice! But maybe we should add why  $w,w',z \notin S \cup \{v,v'\}$? It took me a few minutes to realise why.}.  
 \\
 We observe that there are no edges between $S'$ and $T'$, since $w,w'$ are the only blue vertices in $G[T]$ with a white neighbor. %\julien{I do not understand the use of the last sentence?}. 
 Suppose, towards a contradiction, that $|T'|\leq |S'|$. Then $(T', S')$ forms a valid pair for $\{w,w'\}$. Since the cornerstone $v,v'$ maximises the function $h$, we must have $|T'|-d_{T'} \leq |S|-d_S$. But this is a contradiction, as $S \setminus N[\{v,v',y\}] \subsetneq T' \setminus N[\{w,w',z\}]$.  
\end{proof}
Thus, $(S', T')$ forms a valid pair for the cornerstone $w,w'$, and since the cornerstone $v,v'$ maximises the function $h$, we have $|S'|-d_{S'} \leq |S|-d_S$. Let $\tau_2^{a}$ denote the time needed for $z$ to turn blue, and let $\tau_2^{b}$ denote the time needed thereafter for $|S'|-|S|+d_S \leq d_{S'}$ vertices in $S' \cap N[\{w,w',z\}] $ to be forced blue. Clearly, $\mathbb{E}[\tau_2^{a}|S_{\tau_1},B']=O(1)$. Moreover, by iteratively applying \Cref{cor: color d neighbours in log d steps} to each of the stars centered on $a \in \{w, w', z \}$, with $d=\min\{ |N(a) \cap S'|, |S'|-|S|+d_S \}$, we get $\mathbb{E}[\tau_2^{b}|S_{\tau_1},B']=O(\log(|S'|-|S|+d_S))+O(1)$. 
% \mehdi{We should explain why it is applicable and consider the process restricted to a star. Also below, it is not clear why $\mathbb{E}[\tau_2^{b}|B']=O(\log(|S'|-|S|+d_S))+O(1)$.
% For me $\mathbb{E}[\tau_2^{b}|B']=O(\log(\deg(w)))+O(\log(\deg(w')))+O(\log(\deg(z)))=O(\log(d_{S'}))$ as I apply Corollary 3.10 to the three stars which are restrictions of $N[w],N[w'],N[z]$ once $z$ is colored blue.}
% \julien{This is the thing: we want to turn blue only $|S'|-|S|+d_S$ vertices, so that is why we use Corollary 3.10. If we use the bound that you wrote, we get an additional term that is more than $O(1)$.}\mehdi{to apply Corollary 3.10. you need a star, but then which star do you consider to color only $|S'|-|S|+d_S$ vertices? }
Observe that $\tau_2 \leq \tau_2^{a}+\tau_2^{b}$. Thus, 
\begin{align*}
   \mathbb{E}[\tau_2|S_{\tau_1},B'] \leq\mathbb{E}[\tau_2^{a}|S_{\tau_1},B']+\mathbb{E}[\tau_2^{b}|S_{\tau_1},B'] = O(\log(|S'|-|S|+d_S))+O(1).
\end{align*}
Similarly, we can condition on other possible outcomes $B'' \subseteq B$ and show that $\mathbb{E}[\tau_2|S_{\tau_1},B''] = O(\log(|S'|-|S|+d_S)) +O(1)$. Thus, we have
\begin{align}
\label{eq: boundtau2}
    \mathbb{E}[\tau_2] = \mathbb{E}[\tau_2 \mathbf{1}_{A}]+ \mathbb{E}[\tau_2 \mathbf{1}_{B}]\leq \mathbb{E}[O(\log(|S'|-|S|+d_S))\mathbf{1}_{B}]+O(1).
\end{align}
% \mehdi{Above, I think we rather have $ \mathbb{E}[\tau_2]\le \mathbb{E}[\tau_2 \mathbf{1}_{A}]+ \mathbb{E}[\tau_2 \mathbf{1}_{B}]$ and not necessarily equality as $A$ and $B$ can both occur.}
%In the second sub-phase, we run probabilistic zero forcing on $G[T]$ until $z$ is also blue, which takes $O(1)$ steps in expectation. In the third and final sub-phase, we run probabilistic zero forcing on $G[T]$ until at least $|S_2|-|S|+d_S$ white vertices in $N[\{w,w',z\}]$ are forced blue. We can do this as $d_{S_2} \geq |S_2|-|S|+d_S$. By Corollary \ref{cor: color d neighbours in log d steps}, this takes $O(\log(|S_2|-|S|+d_S))$ steps in expectation. 
%Putting the three sub-phases together, in this case, Step 5 takes at most
%\begin{align*}
    %&\frac{|T|-|S_2|-d_T}{2} + O(\log(|S_2|-|S|+d_S))+ O(1) \\
    %&= \frac{|T|-|S|-d_T+d_S}{2}+O(\log(|S_2|-|S|+d_S))-\left( \frac{|S_2|-|S|+d_S}{2}\right)+ O(1) \\
    %&= \frac{|T|-|S|-d_T+d_S}{2} +O(1)
%\end{align*}
%steps in expectation, as desired.
From \eqref{eq: bound tau1}, we get
\begin{align*}
    \mathbb{E}[\tau_1] &\leq \frac{\mathbb{E}[\min(X_{\tau_1},X_{\tau_1-1}+25)\mathbf{1}_A]+\mathbb{E}[\min(X_{\tau_1},X_{\tau_1-1}+25)\mathbf{1}_B]}{2}-\frac{d_T}{2}\\
    & \leq \frac{\mathbb{E}[(X_{\tau_1-1}+25)\mathbf{1}_A]+\mathbb{E}[X_{\tau_1}\mathbf{1}_B]}{2}-\frac{d_T}{2}\leq \frac{\mathbb{E}[X_{\tau_1-1}\mathbf{1}_A]+\mathbb{E}[X_{\tau_1}\mathbf{1}_B]}{2}-\frac{d_T}{2}+O(1) \\
    & \leq \frac{\mathbb{E}[(|T|-|S|+d_S)\mathbf{1}_A]+\mathbb{E}[(|T|-|S'|)\mathbf{1}_B]}{2}-\frac{d_T}{2}+O(1) \\
    &= \frac{|T|-d_T-|S|+d_S}{2}-\mathbb{E}\left[\left(\frac{|S'|-|S|+d_S}{2}\right)\mathbf{1}_B\right]+O(1).
\end{align*}
% \mehdi{Why does the last equality hold? Also, for $\mathcal{G}$ a sigma-algebra, $\bE[\tau_2|\mathcal{G}]$ is the unique rv satisfying $\bE[\tau_21_{B}]=\bE[\bE[\tau_2|\mathcal{G}]1_{B}]$. Which sigma-algebra $\mathcal{G}$ do you take? }
% \julien{There is no need to reason with sigma-algebras here: the last equality is simply by plugging $\mathbf{1}_A+ \mathbf{1}_B=1$. But if you really want to reason with sigma-algebras, this is simply the sigma-algebras generated by the configurations of blue vertices at each step of the process.}
Combining with \eqref{eq: boundtau2}, we obtain 
\begin{align*}
    \mathbb{E}[Z_5] &= \mathbb{E}[\tau_1]+\mathbb{E}[\tau_2] \\
    &\leq \frac{|T|-d_T-|S|+d_S}{2}+ \mathbb{E}\left[\left(O(\log(|S'|-|S|+d_S))-\frac{|S'|-|S|+d_S}{2}\right)\mathbf{1}_B\right]+O(1).
\end{align*}
Since for every absolute constant $\eta > 0$, there exists an absolute constant $\gamma$ such that $\eta \log x -x/2 \leq \gamma$ for every $x >0$, we get
\begin{align*}
    \mathbb{E}[Z_5] & \leq \frac{|T|-d_T-|S|+d_S}{2}+O(1).
\end{align*}
% \mehdi{why does the last inequality hold? Note that $S'$ is a random variable.}
% \julien{I have added a sentence to explain}
This proves the lemma.
\end{proof}

Note that at least $d_S$ vertices in $S$ were colored blue during Step 4, and therefore both $G[S]$ and $G[T]$ have at most $|S|-d_S$ white vertices at the beginning of Step 6. The proof of the next lemma then follows immediately from the proof of \cite[Lemma 4.10]{Sun}.
%\julien{Say the proof of the next lemma comes immediately from \cite{Sun}.}

% \begin{lemma}
%     \label{lem:Bounds-Step6}
%     Step 6 takes at most $\mathbb{E}[\log (|S_2|-(|S|-d_S))]+O(1)$ steps in expectation.
% \end{lemma}

\begin{lemma}
    \label{lem:Bounds-Step6}
    Step 6 takes at most $|S|-d_S+O(1)$ steps in expectation.
\end{lemma}

It is then easy to deduce Theorem~\ref{ept improved} from Lemmas \ref{lem:Bounds-Step4}, \ref{lem:Bounds-Step5} and \ref{lem:Bounds-Step6}.
\begin{proof}[Proof of Theorem~\ref{ept improved}]
Let $Z_4$, $Z_5$, $Z_6$ be the random variables which denote the number of steps for which Step 4, 5 and 6 respectively last, and let $\Gamma$ be the total time duration of probabilistic zero forcing on $G$.
Then by applying \Cref{lem:Bounds-Step4,lem:Bounds-Step5,lem:Bounds-Step6}, we obtain
\begin{align*}
    \mathbb{E}[\Gamma] &\leq \mathbb{E}[Z_4]+\mathbb{E}[Z_5]+\mathbb{E}[Z_6] \\
    &\leq O(\log (d_T+d_S))+O(1) + \frac{|T|-d_T-|S|+d_S}{2}+O(1) + |S|-d_S+O(1) \\
    &= \frac{|T|+|S|}{2} + O(\log (d_T+d_S))-\left(\frac{d_T+d_S}{2} \right) + O(1) \\
    &= \frac{n}{2}+O(1),
\end{align*}
as wanted.
\end{proof}

%\nocite{*}
\bibliographystyle{abbrv}
\renewcommand{\bibname}{Bibliography}
\bibliography{bib}

\appendix
\section{Proofs of Lemma \ref{lem: coloring d neighbours in a star} and Corollary \ref{cor: color d neighbours in log d steps}}
\label{Appendix-Bounds-Star}

%\begin{proof}[Proof of \Cref{exp decay to color N(v)}]
%First of all, by \Cref{lemma 2.10 for Sun}, we may assume that initially only $v$ is coloured blue.
%We then consider the modified probabilistic process where only vertex $v$ is allowed to propagate its colour to other vertices, with same probability as in the probabilistic zero forcing.
%This modified probabilistic process trivially satisfies the assumptions of \Cref{lemma 2.11 for Sun}, and therefore it suffices to show the tail bound for this modified probabilistic process to establish the first part of the lemma.
%Note that this modified probabilistic process is simply a probabilistic zero forcing propagation on a star with $d$ leaves, starting with the centre vertex coloured blue.
%\Cref{color the star} then provides the desired bound. \\
%For the second part, we have 
%\begin{align*}
    %\mathbb{E}[\Gamma 1_{\Gamma \geq C \log d}] \leq \sum_{t \geq C \log d} (t+1) \mathbb{P}(\Gamma \geq t) = O\left (\sum_{t \geq C \log d} t \alpha^t\right) =o(1),
%\end{align*}
%and therefore
%\begin{align*}
    %\mathbb{E}[\Gamma] = \mathbb{E}[\Gamma 1_{\Gamma < C \log d}] +\mathbb{E}[\Gamma 1_{\Gamma \geq C \log d}] = O(\log d).
%\end{align*}
%\end{proof}
\begin{proof}[Proof of Lemma \ref{lem: coloring d neighbours in a star}]
    We partition the interval $[0,n]$ into subintervals as follows. Let $I_1 =[0,1)$. If $I_j=[a_j,b_j)$ for $b_j<n/3$, we set $I_{j+1}=\left[b_j, b_j +\frac{b_{j}+1}{6}\right) \cap \left[0,\frac{n}{3}\right)$. If $I_J=[a_J,b_J)$ for $b_J=\frac{n}{3}$, then set $I_{J+r}=\left[n-\frac{2n}{3}\left(\frac{5}{6}\right)^{r-1},n-\frac{2n}{3}\left(\frac{5}{6}\right)^{r}\right)$. For the least value of $R$ such that $\frac{2n}{3}\left( \frac{5}{6}\right)^{R-1}<1$, we set $I_{J+R}=\left[n-\frac{2n}{3}\left( \frac{5}{6}\right)^{R-1},n\right]$ to be the final interval.
    Note that $I_1,\dots,I_J$ partition $\left[0,\frac{n}{3}\right)$, and $I_{J+1},\dots,I_{J+R}$ partition $\left[\frac{n}{3},n\right]$, so we have a complete partition of $[0,n]$. \\
    It is easy to see that $J \leq C_1 \log n$ and $R \leq C_2 \log n$ for some $C_1, C_2>0$.
    If $d \leq n-1$, let $m$ be the smallest integer such that the left end point of $I_m$, $a_m$ satisfies $a_m \geq d$. Otherwise, let $m=J+R$. If $d \leq \frac{n}{3}$, it is easy to see that $m \leq C_3 \log d$, for some constant $C_3>0$. On the other hand, if $d \geq \frac{n}{3}$, then $m \leq J+R \leq (C_1+C_2)\log n \leq (C_1+C_2)\log (3d)\leq C_4 \log d$, for some $C_4>0$. So, we have $m \leq C'\log d$, where $C'=\max(C_3,C_4)$.\\ Note that, by \cite[Lemma 3.2]{Sun}, if the number of blue leaves is $k \in I_r$ for $r \leq J$, with probability at least $1/5$, the number of blue leaves after the next iteration of probabilistic zero forcing will be in some $I_s$ for $s>r$. Moreover, by \cite[Lemma 3.3]{Sun}, the same is true for $J+1 \leq r < J+R$. \\
    Therefore, provided $t \geq 10C'\log d$, the probability that the blue vertex will propagate to fewer than $d$ leaves in $t$ steps is at most
    \begin{align*}
    \mathbb{P}\left(\text{Bin}\left(t,\frac{1}{5}\right) \leq m\right) \leq \mathbb{P}\left(\text{Bin}\left(t,\frac{1}{5}\right) \leq C'\log d\right) &\leq \mathbb{P}\left(\text{Bin}\left(t,\frac{1}{5}\right) \leq \frac{t}{10}\right) \leq e^{-\frac{t}{40}},     
    \end{align*}
    where the last inequality follows from \Cref{thm:Chernoff}. Thus, if we set $\alpha=e^{-1/40}$ and $C=10C'$, then the probability that at least $d$ leaves are blue after $t>C\log d$ steps is at least $1-\alpha^t$, as claimed.
\end{proof}

\begin{proof}[Proof of Corollary \ref{cor: color d neighbours in log d steps}]
% \mehdi{remark: I think we no longer need to prove this result. I found that it is proved in Theorem 2.7, in the paper "Propagation time for probabilistic zero forcing" written by Geneson and Hogben. But the proof is short, we can keep it if you want.}
% \julien{I do not think it is correct: the theorem you refer to is the expected time for turning the full star blue, whereas we only want part of it here, and the bound for the whole star is too crude fot our application.}
    Using the tail sum formula for expectation, we have
    \begin{align*}
        \mathbb{E}[\Gamma]=\sum_{t=1}^{\infty}\mathbb{P}(\Gamma \geq t)&= \sum_{t \leq C \log d} \mathbb{P}(\Gamma \geq t)+\sum_{t > C\log d}\mathbb{P}(\Gamma \geq t)\\
        &\leq \sum_{t \leq C \log d} \mathbb{P}(\Gamma \geq t)+\sum_{t > C\log d} \alpha^t = O(\log d)+O(1)=O(\log d),
    \end{align*}
    where the first inequality follows from Lemma \ref{lem: coloring d neighbours in a star}.
\end{proof}

\section{Proof of Lemma \ref{lem:Bounds-Step4}}
\label{Appendix-Bounds-Step4}

  \begin{proof}[Proof of \Cref{lem:Bounds-Step4}]
  In this proof, we only treat the case where our choice in Step 1 was a cornerstone $v,v'$ with $y$ as the choice of the common neighbor: the other cases can be deduced similarly via a slightly simpler proof.
    Let $\Lambda$ be the random variable associated to the number of rounds that Step 4 lasts.
    We consider the following stopping times: let $\tau_1$ be the first time at which all of $N[v]$ is colored blue, let $\tau_2$ be the first time after $\tau_1$ (i.e. $\tau_2 \geq \tau_1$) at which all of $N[y]$ is colored blue, and let $\tau_3$ be the first time after $\tau_2$ (i.e. $\tau_3 \geq \tau_2$) at which all of $N[v']$ is colored blue.
    Let $\Lambda_1=\tau_1$, $\Lambda_2=\tau_2-\tau_1$, and $\Lambda_3=\tau_3-\tau_2$.
    By \Cref{lemma 2.11 for Sun}, we have that $\mathbb{E}[\Lambda_1]$ is at most $\ept(S_{d(v)},\{w\})$  %\julien{Make sure this definition is understood}
    , where $d(v)$ is the degree of $v$ and $S_{d(v)}$ is a star on $d(v)+1$ vertices centered on a vertex $w$.
    By Corollary \ref{cor: color d neighbours in log d steps}, we obtain
    \begin{align*}
        \mathbb{E}[\Lambda_1]=O(\log d(v)) +O(1)= O(\log (d_S+d_T)) +O(1), 
    \end{align*}
    where the term $O(1)$ accounts for the case where $d(v)=1$.
    % \mehdi{why is the latter needed?} \julien{Because then $\log 1=0$ but the expected runtime is $1$}
    We obtain similar bounds for $\mathbb{E}[\Lambda_2]$ and $\mathbb{E}[\Lambda_3]$, from which we conclude that
    \begin{align*}
        \mathbb{E}[\Lambda]=\mathbb{E}[\Lambda_1]+\mathbb{E}[\Lambda_2]+\mathbb{E}[\Lambda_3]=O(\log (d_S+d_T)) +O(1).
    \end{align*}
\end{proof}

\end{document}